\documentclass[10pt]{article}
\usepackage[utf8]{inputenc}
\usepackage[T2A]{fontenc}
\usepackage[english]{babel}
\usepackage{amsmath,amssymb,amsthm}
\usepackage{graphicx}
\usepackage{pdflscape}
\usepackage{url}
\usepackage{hhline}

\title{Surfaces containing two isotropic circles through each point}
\author{Egor Morozov\thanks{National Research University Higher School of Economics, Faculty of Mathematics, Russian Federation\endgraf\textit{E-mail:} \texttt{eamorozov\_1@edu.hse.ru}}}
\date{}

\DeclareMathOperator{\const}{const}
\DeclareMathOperator{\id}{Id}
\DeclareMathOperator{\cl}{Cl}

\newcommand{\R}{\mathbb R}
\newcommand{\mattc}{\mathrm{Mat}_2(\mathbb C)}
\newcommand{\gltc}{\mathrm{GL}_2(\mathbb C)}
\newcommand*{\hm}[1]{#1\nobreak\discretionary{}%
{\hbox{$\mathsurround=0pt #1$}}{}}
\renewcommand{\C}{\mathbb C}

\theoremstyle{plain}
\newtheorem{lemma}{Lemma}
\newtheorem{proposition}{Proposition}
\newtheorem{corollary}{Corollary}
\newtheorem{theorem}{Theorem}
\newtheorem{conjecture}{Conjecture}
\newtheorem{problem}{Problem}

\theoremstyle{remark}
\newtheorem{remark}{Remark}
\newtheorem{example}{Example}

\begin{document}
\maketitle
\begin{abstract}
We prove (under some technical assumptions) that each surface in $\mathbb R^3$ containing two arcs of parabolas with axes parallel to $Oz$ through each point has a parametrization $\left(\frac{P(u,v)}{R(u,v)},\frac{Q(u,v)}{R(u,v)},\frac{Z(u,v)}{R^2(u,v)}\right)$ for some $P,Q,R,Z\in\mathbb R[u,v]$ such that $P,Q,R$ have degree at most 1 in $u$ and $v$, and $Z$ has degree at most 2 in~$u$ and~$v$. The proof is based on the observation that one can consider a parabola with vertical axis as an isotropic circle; this allows us to use methods of the recent work by M.\,Skopenkov and R.\,Krasauskas in which all surfaces containing two Euclidean circles through each point are classified. Such approach also allows us to find a similar parametrization for surfaces in $\mathbb R^3$ containing two arbitrary isotropic circles through each point (under the same technical assumptions). Finally, we get some results concerning the top view (the projection along the $Oz$ axis) of the surfaces in question.

\smallskip


\smallskip

\end{abstract}
 
\section{Introduction and background}
Various problems about classification of surfaces in $\R^3$ containing two special curves through each point are classical and have a long history. Probably the most classical example is a one-sheeted hyperboloid containing two lines through each point (this fact was discovered by C.\,Wren back in 1669 and later used by V.\,Shukhov in his famous hyperboloid structures). Another example is the problem of finding all surfaces containing several circles through each point. Although the problem has been studied since XIXth century, it was solved only in the recent paper~\cite{skopenkov2019surfaces} (a~brief but informative history of the question can be also found there).

\begin{figure}[h]
\begin{minipage}{0.32\linewidth}
\center{\includegraphics[width=1\linewidth]{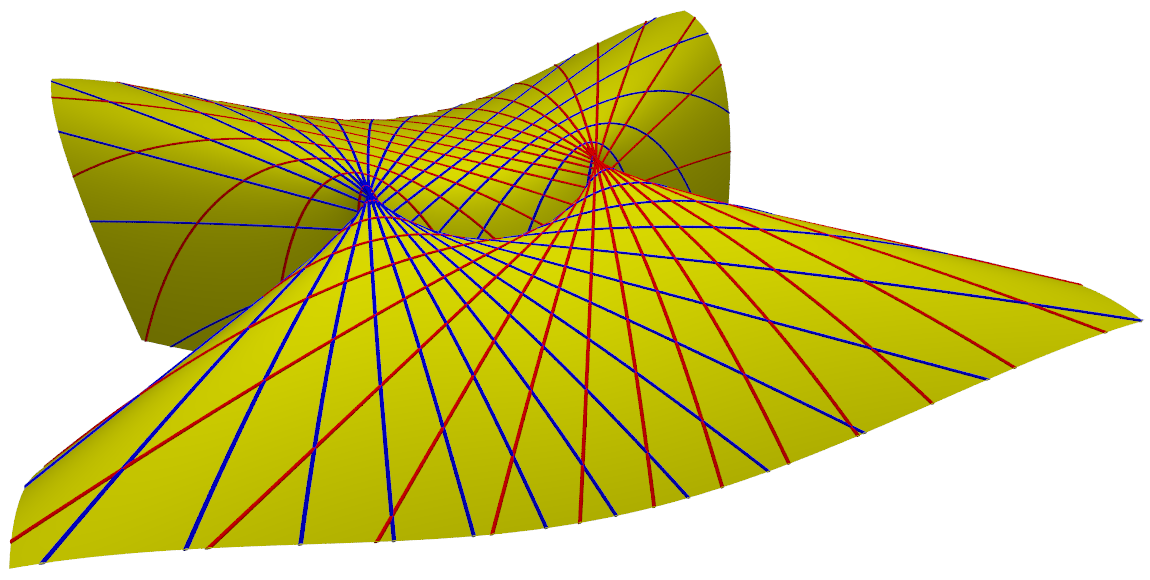}} \\ $a$
\end{minipage}
\hfill
\begin{minipage}{0.32\linewidth}
\center{\includegraphics[width=1\linewidth]{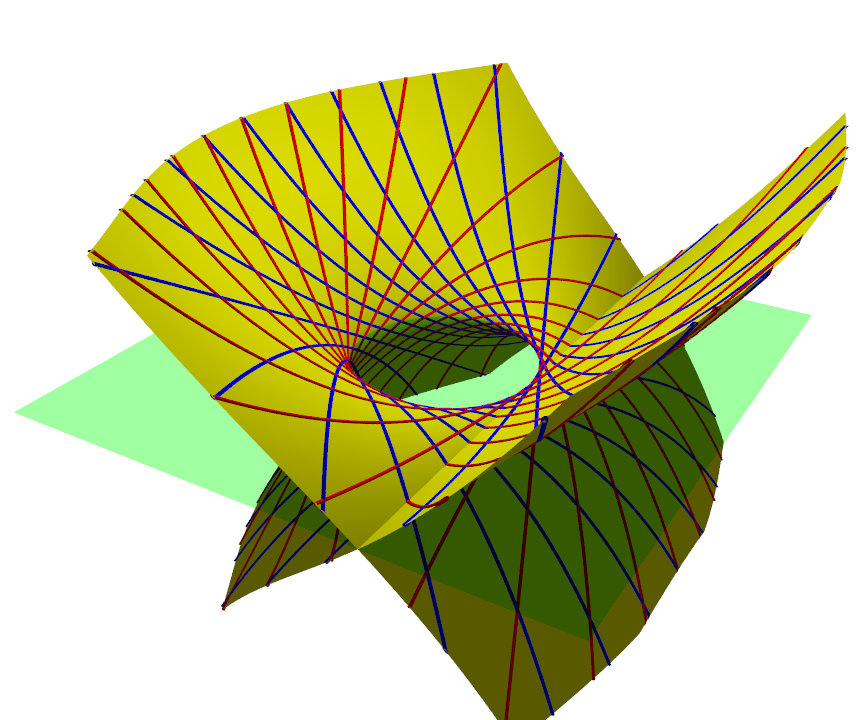}} \\ $b$
\end{minipage}
\hfill
\begin{minipage}{0.32\linewidth}
\center{\includegraphics[width=1\linewidth]{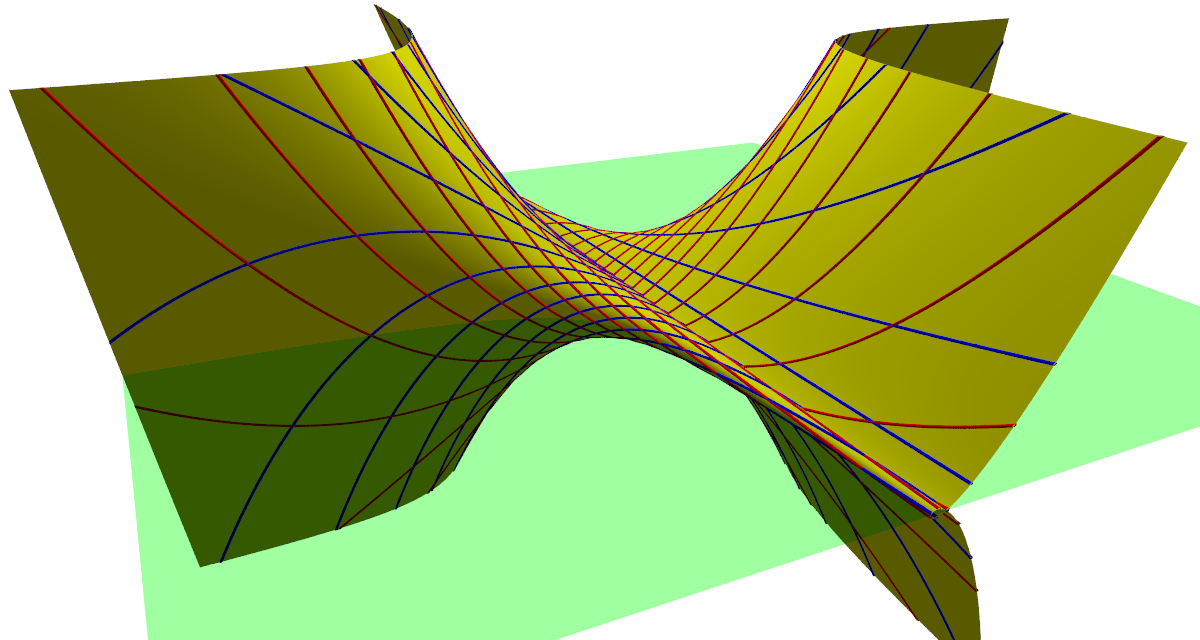}} \\ $c$
\end{minipage}
\caption{examples of surfaces containing two parabolas with vertical axes through each point. In figures $b,c$ the $xy$-plane is colored in green.}
\label{fig:examples-1}
\end{figure}

\begin{figure}[h]
\begin{minipage}{0.35\linewidth}
\center{\includegraphics[width=1\linewidth]{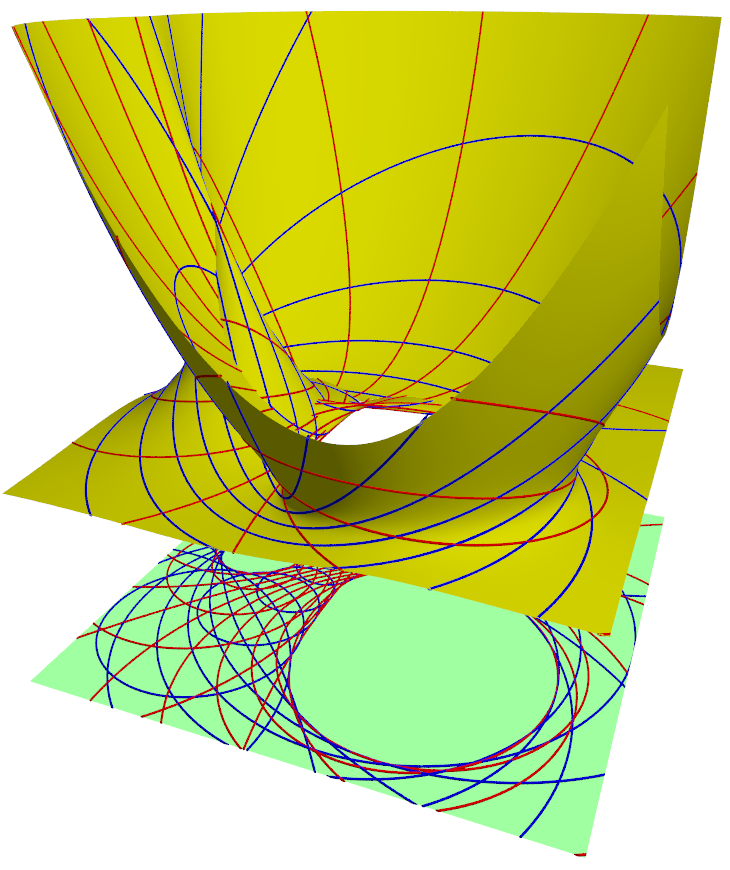}}
\end{minipage}
\hfill
\begin{minipage}{0.45\linewidth}
\center{\includegraphics[width=1\linewidth]{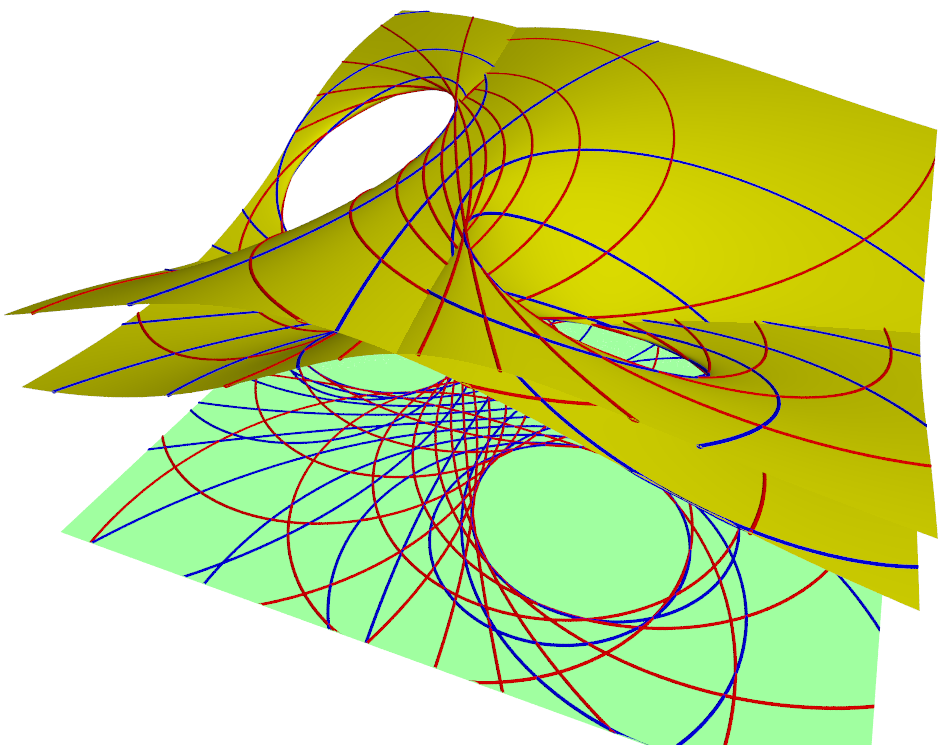}}
\end{minipage}
\caption{examples of surfaces containing two general isotropic circles through each point and the top views of isotropic circles.}
\label{fig:examples-2}
\end{figure}

The aim of this paper is finding all surfaces in $\R^3$ containing two parabolas with vertical axes through each point (Theorem~\ref{thm:main1}; see Fig.~\ref{fig:examples-1} for some examples). The new idea is to consider parabolas with vertical axes as circles in isotropic geometry and apply isotropic stereographic projection (see~\S\S\ref{sec:statements},\ref{sec:outline} and~\cite{krasauskas2014bilinear,pottmann2007discrete}). This approach allows us to solve the problem using quite standard methods similar to~\cite{skopenkov2019surfaces}. In passing, we find all surfaces containing two arbitrary isotropic circles through each point (Theorem~\ref{thm:main2}; see Fig.~\ref{fig:examples-2} for some examples). Particular cases called \emph{isotropic cyclides} are considered in~\cite{dahl2014improved,krasauskas2014bilinear,sachs1990isotrope}.

An interesting feature of isotropic geometry is the concept of ``top view''. Projecting the isotropic circles through each point of the surface to a horizontal plane, we obtain a configuration of lines and (Euclidean) circles in the plane (see Fig.~\ref{fig:examples-2}). For the case of surfaces containing two parabolas through each point we give a complete classification of the resulting line configurations (Corollary~\ref{cor:top1}) and for the case of surfaces containing two isotropic circles through each point we give a partial classification (Corollary~\ref{cor:top2}). It seems to be not difficult to obtain a complete classification for the latter case using the methods of \S\S\ref{sec:linal}---\ref{sec:top2} but this is too technical and lies beyond the scope of the paper (see Conjecture~\ref{con:topview} in \S\ref{sec:open}).

Our proofs use Schicho's theorem (Theorem~\ref{thm:schicho}) which parametrizes all surfaces containing two conics through each point. Although Schicho's theorem is powerful, it requires many technical assumptions (e.~g., the surface contains finitely many conics through each point) which leads to quite long statements of our main theorems. We give a detailed plan of the proof in \S\ref{sec:outline}.

The problem of finding all surfaces containing two circles through each point is interesting also in higher dimensions and other geometries (for instance, classification of surfaces in $\R^3$ was actually deduced from the classification of surfaces in $\R^4$ in~\cite{skopenkov2019surfaces}). A closely related problem is finding all surfaces containing a line and a circle through each point. By ``geometry'' we mean only how the notions of ``circle'' and ``line'' are defined; in particular, we do \emph{not} study classification of surfaces up to a transformation group; we are interested just in finding all the surfaces in question.

\begin{landscape}
\begin{table}[h]
\caption{a summary of known results on surfaces containing two special curves through each point in various dimensions and geometries.}
\label{tab:all-results}
\begin{center}
\begin{tabular}{|p{1.2cm}|p{4cm}|p{6cm}|p{2.8cm}|p{3.5cm}|}
\hline
Problem & Assumptions & Readily-calculable solution & Author(s) & Reference \\
\hline
CCE3 & --- & yes & Skopenkov-Krasauskas & \cite[Theorem~1.1]{skopenkov2019surfaces} \\
\hhline{|-|-|-|~|-|}
CCE4 & two circles through each point are noncospheric, finitely many circles through each point & no & & \cite[Corollary~5.1]{skopenkov2019surfaces} \\
\hhline{|~|-|-|-|-|}
 & infinitely many circles through each point & yes & Koll\'ar & \cite[Remark~5.3]{skopenkov2019surfaces} summarizing \cite[Theorems~3, 8 and Propositions~11, 12]{kollar2018quadratic} \\
\hhline{|~|-|-|-|-|}
 & the surface is algebraic & no & Lubbes & \cite[Theorem~1]{lubbes2021surfaces} \\
\hline
CCE9, CCI9 & --- & yes (there are no such surfaces, i.~e., any surface in question is contained in $\R^8$) & Schicho & remark after \cite[Theorem~4.1]{skopenkov2019surfaces} \\
\hline
CCI3 & two isotropic circles through each point are noncospheric, finitely many isotropic circles through each point & yes & \multicolumn{2}{c|}{see Theorem~\ref{thm:main2} below} \\
\hline
CLE3 & --- & yes (such surfaces are quadrics only) & Nilov-Skopenkov & \cite[Theorem~1.1]{nilov2013surface} \\
\hline
CLE4, CLE5 & the surface is algebraic & no & Lubbes & \cite[Corollary~4]{lubbes2021surfaces} \\
\hhline{|-|~|-|~|~|}
CLE6 & & yes (there are no such surfaces, i.~e., any surface in question is contained in $\R^5$) & & \\
\hhline{|-|-|-|~|-|}
CLS3 & the surface is algebraic and contains \emph{exactly} one circle and one line through each point & no & & \cite[Corollary~3b]{lubbes2019euclidean} \\
\hline
\end{tabular}
\end{center}
\end{table}
\end{landscape}

State of the art on both problems is summarized in Table~\ref{tab:all-results}. To make the table visible, we have to focus just on the very two particular problems and also to omit many remarkable results covered by the ones presented in the table. We refer to \cite{lubbes2021surfaces} for a more complete survey and a discussion of very interesting closely related questions.

\textbf{How to read Table~\ref{tab:all-results}.} Each row depicts a problem of finding all surfaces containing two special curves through each point in space of given dimension with given geometry. The columns mean the following:

\begin{itemize}
\item{\textit{Problem.} A label in this column abbreviates a particular problem as follows.
\begin{itemize}
\item{The first two letters refer to curves passing through each point of the surface: ``CC'' means ``two circles'' and ``CL'' means ``a line and a circle''.}
\item{The third letter refers to how the notions of ``circle'' and ``line'' are understood: ``E'' means Euclidean circle and line in $\R^n$, ``I'' means isotropic circle and line in $\R^n$ (see the beginning of \S\ref{sec:statements} and~\cite{sachs1990isotrope} for the case $n=3$), ``S'' means circle and spherical line (i.~e. great circle) in $S^n$. Note that the problems ``CCEn'' and ``CCSn'' are equivalent, but this is not the case for ``CLEn'' and ``CLSn''. Similar problems in pseudoeuclidean and hyperbolic geometries are unsolved.}
\item{The number refers to the dimension of the ambient space; we assume that the surface is not contained in a subspace of lower dimension.}
\end{itemize}
For example, the label ``CCE3'' denotes a problem of finding all surfaces in $\R^3$ which are not contained in $\R^2$ and which contain two Euclidean circles through each point.}
\item{\textit{Assumptions.} These are the assumptions under which the solution has been obtained. The analyticity of the surface and the analytical dependence of the curves through each point on the point are the default ones. It seems that analyticity is not really a restriction (see the paragraph before~\cite[Problem~5.7]{skopenkov2019surfaces}). We assume it in order to avoid technical difficulties. At the same time, we believe that the analyticity assumptions are general enough to make results interesting and useful for applications.}
\item{\textit{Readily-calculable solution.} Informally, we say that a problem has a \emph{readily-calculable solution} if there is an explicit \emph{parametrization of the set} of surfaces in question, i. e. an explicit surjection from a $p$-dimensional subset in $\R^p$ (or a finite union of such subsets with different $p$) to the set of surfaces in question. This is similar to solution of an equation: if the number of solutions is finite, then one can just list them; otherwise one parametrizes the set of solutions. This notion pretends to be neither precise nor ideal nor universal. ``Yes'' standing in the 2nd column of Table 1 is just an invitation for a \emph{user} to apply the cited solution immediately, and ``no'' is an invitation for a \emph{developer} to make his or her own research on the particular problem. Results marked by ``yes'' give an impression of ``the best possible answers'' to the problems. But this does not pretend to be the case outside the context of Table 1 (and we avoid further specification of the notion of a ``readily-calculable solution'' because it would not change Table 1). See \cite[Remark~1.2]{skopenkov-embeddings} for a detailed discussion.}
\end{itemize}

The problems in Table~\ref{tab:unsolved} have not yet been solved in any sense but seem to be accessible by known methods.

\begin{table}[h]
\caption{Remarks on some unsolved problems.}
\label{tab:unsolved}
\begin{center}
\begin{tabular}{|p{1.4cm}|p{11cm}|}
\hline
Problem & Remark \\
\hline
CCE8 & The problem might possibly be solved via octonions (just as CCE4 was approached by quaternions, see \cite{skopenkov2019surfaces}).\\
\hline
CCI5 & Using the method of the paper the problem can be reduced to the classification of Pythagorean 6-tuples from~\cite[Theorem~1.3]{skopenkov2019surfaces}.\\
\hline
\end{tabular}
\end{center}
\end{table}

The 3-dimensional pseudoeuclidean case may be interesting due to relation to \emph{circle patterns} (see~\cite{pottmann2007discrete} for the isotropic case).

\textbf{Applications.} The problems in question are shown to be useful in various applications. The book~\cite{pottmann2001computational} contains detailed study of ruled surfaces modeling. Classical examples of surfaces containing several Euclidean circles through each point are cyclides, which have numerous applications in architecture~\cite{pottmann2012darboux}, CAGD (see references in~\cite{krasauskas2000studying}) and in the study of the kinematic map~\cite{krasauskas2020kinematic,lubbes2018kinematic}. For applications of isotropic cyclides in CAGD see~\cite{dahl2014improved,krasauskas2014bilinear}. General surfaces containing isotropic circles through each point are connected with the study of surfaces enveloped by 1-parametric families of rotational cones. This connection is made through Laguerre geometry and finds unexpected applications in computer numerically controlled machining~\cite{skopenkov2020characterizing} (see also Problem~\ref{prb:laguerre} in \S\ref{sec:open}). Note that isotropic model of Laguerre geometry is also applied to the study of Dupin cyclides~\cite{krasauskas2000studying,pottmann1998applications}.

\textbf{Organization of the paper.} In \S\ref{sec:statements} we state the main results and in \S\ref{sec:outline} we state a number of auxiliary ones proved in \S\S\ref{sec:goodparam}-\ref{sec:polynomials}. Theorems~\ref{thm:main1} and \ref{thm:main2} are proved in \S\S\ref{sec:main1}-\ref{sec:main2} respectively. Corollary~\ref{cor:top1} is proved in \S\ref{sec:top1}. Corollary~\ref{cor:top2} is much more complicated; it is proved in \S\ref{sec:top2} using a small theory developed in \S\S\ref{sec:linal}-\ref{sec:top2}. Finally, in \S\ref{sec:open} we state several open problems.

\section{Statements}\label{sec:statements}
Consider 3-dimensional Euclidean space $\R^3$ with the Cartesian coordinates $(x,y,z)$. Each line parallel to $Oz$ is called \emph{vertical}. Each plane parallel to the plane $Oxy$ is called \emph{horizontal}.

\begin{figure}[h!]
\center{\includegraphics[scale=0.5]{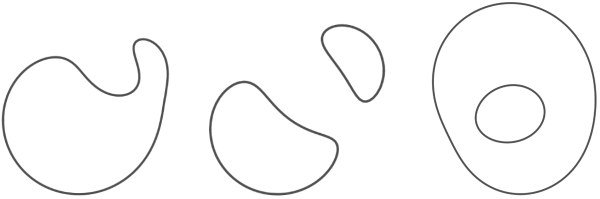}}
\caption{examples of cyclics.}
\label{fig:cyclics}
\end{figure}

An \emph{isotropic circle} in $\R^3$ is either a parabola with vertical axis or an ellipse whose orthogonal projection to the plane $Oxy$ is a circle. The latter projection is called \emph{the top view} of the isotropic circle. An \emph{isotropic sphere} in $\R^3$ is either a paraboloid of revolution with vertical axis or a circular cylinder with vertical axis. 
A \emph{pencil of lines} in the plane is either a set of all lines passing through some fixed point or a set of all lines parallel to some fixed line. A \emph{cyclic} is a subset of the plane given by the equation
\begin{equation}\label{eqn:cyclic-general}
a(x^2+y^2)^2+(x^2+y^2)(bx+cy)+Q(x,y)=0,
\end{equation}
where $a,b,c\in\R$ and $Q\in\R[x,y]$ has degree at most 2 (and $a,b,c,Q$ do not vanish simultaneously). Figure~\ref{fig:cyclics} gives an idea of how cyclics look like.

An \emph{analytic surface} in $\R^3$ is the image of an injective real analytic map of a planar domain into $\R^3$ with nondegenerate differential at each point. An isotropic circular arc \emph{analytically dependent} on a point is a real analytic map of an analytic surface into the real analytic variety of all isotropic circular arcs in $\R^3$.

Denote by $\R_{i,j}$ (respectively, $\C_{i,j}$) the set of all polynomials $F\in\R[u,v]$ (respectively, $\C[u,v]$) such that $\deg_u F\le i$ and $\deg_v F\le j$.

\begin{theorem}\label{thm:main1}
Assume that through each point of an analytic surface in $\R^3$ one can draw two transversal parabolic arcs with vertical axes, fully contained in the surface (and analytically depending on the point). Assume that these two arcs lie neither in the same isotropic sphere nor in the same plane. Assume that through each point in some dense subset of the surface one can draw only finitely many (not nested) arcs of isotropic circles and line segments contained in the surface. Then the surface (possibly besides a one-dimensional subset) has a parametrization
\begin{equation}\label{eqn:mainparam1}
\Phi(u,v)=\left(\frac{P(u,v)}{R(u,v)},\frac{Q(u,v)}{R(u,v)},\frac{Z(u,v)}{R^2(u,v)}\right)
\end{equation}
for some $P,Q,R\in\R_{1,1}$ and $Z\in\R_{2,2}$, where $R\ne0$, such that the parabolic arcs are the curves $u=\const$ and $v=\const$.
\end{theorem}

\begin{figure}[h]
\begin{minipage}{0.3\linewidth}
\center{\includegraphics{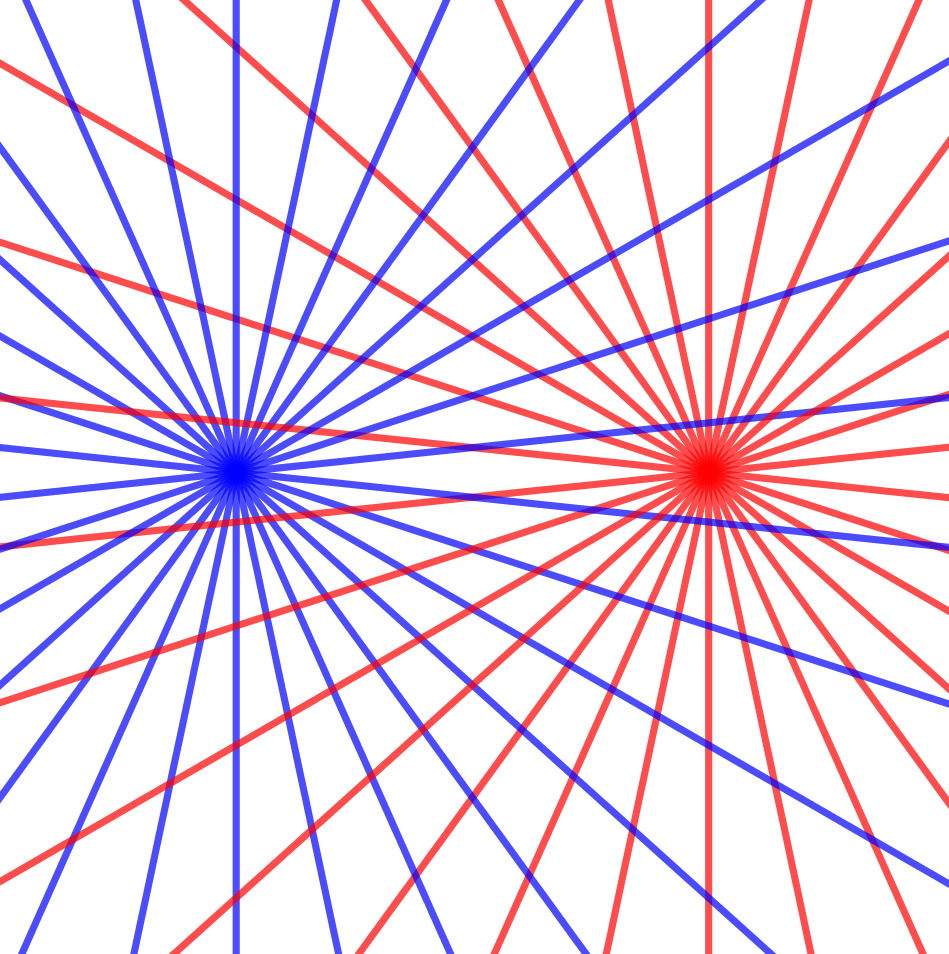}} \\ $a$
\end{minipage}
\hfill
\begin{minipage}{0.3\linewidth}
\center{\includegraphics{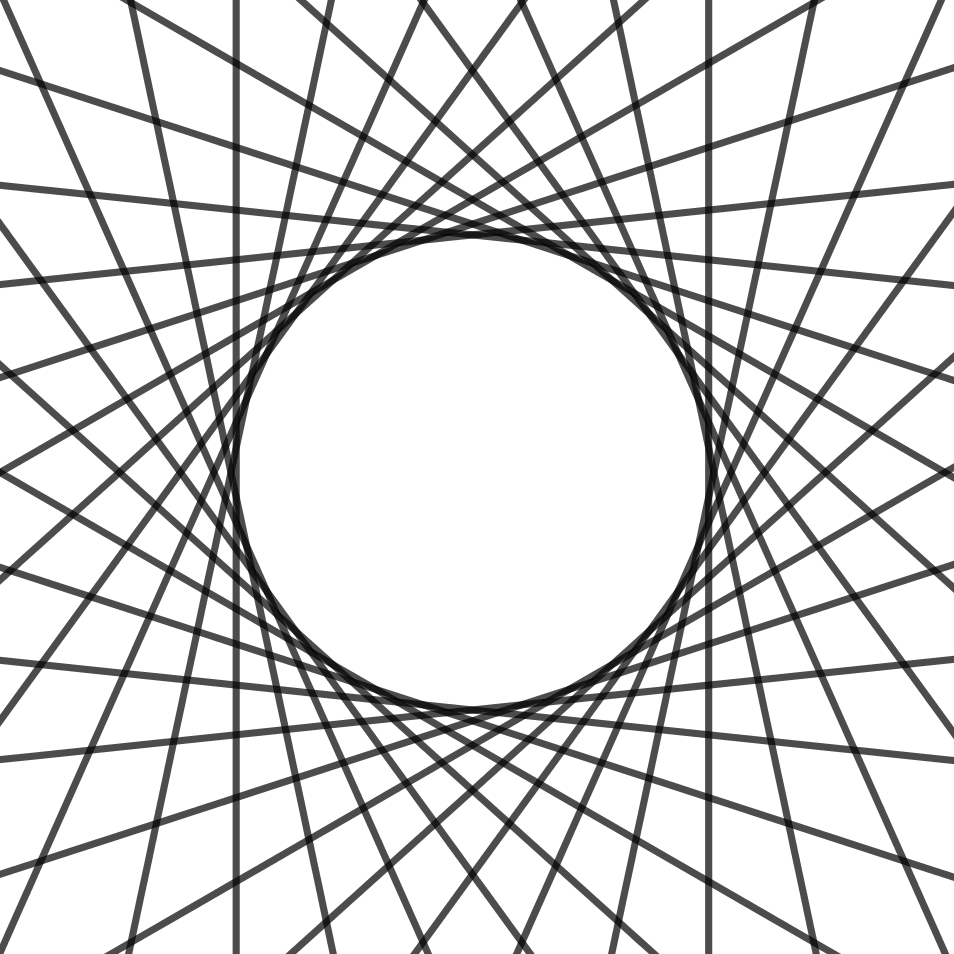}} \\ $b$
\end{minipage}
\hfill
\begin{minipage}{0.3\linewidth}
\center{\includegraphics{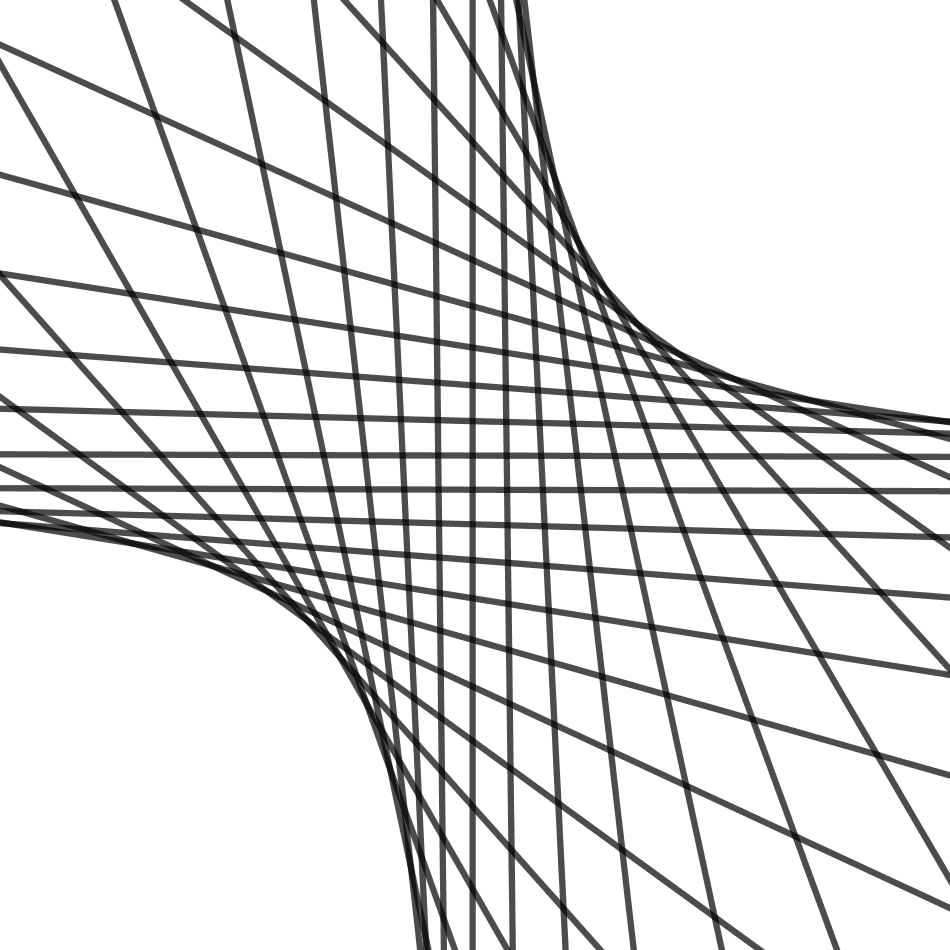}} \\ $c$
\end{minipage}
\caption{to Corollary~\ref{cor:top1}: the top views of parabolas through each point.}
\label{fig:top1}
\end{figure}

\begin{corollary}\label{cor:top1}
For each surface satisfying the assumptions of Theorem~\ref{thm:main1} the top views of the two parabolas through each point are either tangent to one conic or lie in a union of two pencils of lines (see Fig.~\ref{fig:top1}).
\end{corollary}

\begin{remark}
Actually all figures in Fig.~\ref{fig:top1} depict a set of lines dual to the points of some conic. However, in Fig.~\ref{fig:top1}a the conic is not smooth. This is the case for the surface shown in Fig.~\ref{fig:examples-1}a; the base points of both pencils are the projections of the surface singularities. In Fig.~\ref{fig:examples-1}b,c the top views of parabolas through each point are lines tangent to a circle and a hyperbola respectively (see Fig.~\ref{fig:top1}b,c).
\end{remark}

\begin{figure}[h]
\begin{minipage}{0.48\linewidth}
\center{\includegraphics[width=1\linewidth]{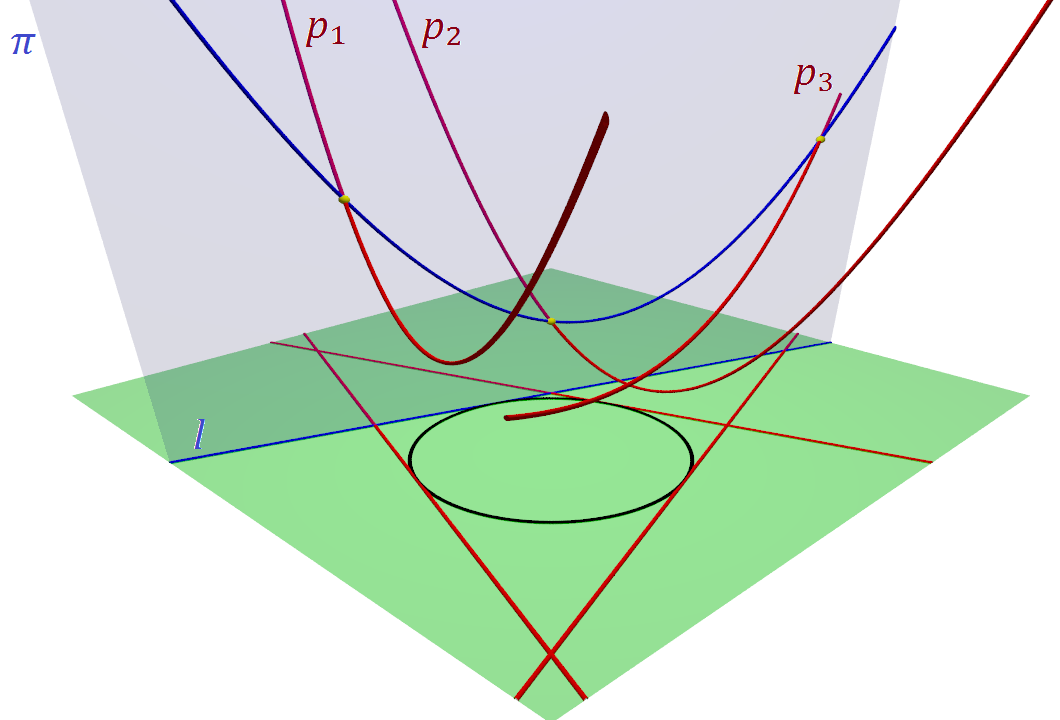}}
\end{minipage}
\hfill
\begin{minipage}{0.48\linewidth}
\center{\includegraphics[width=1\linewidth]{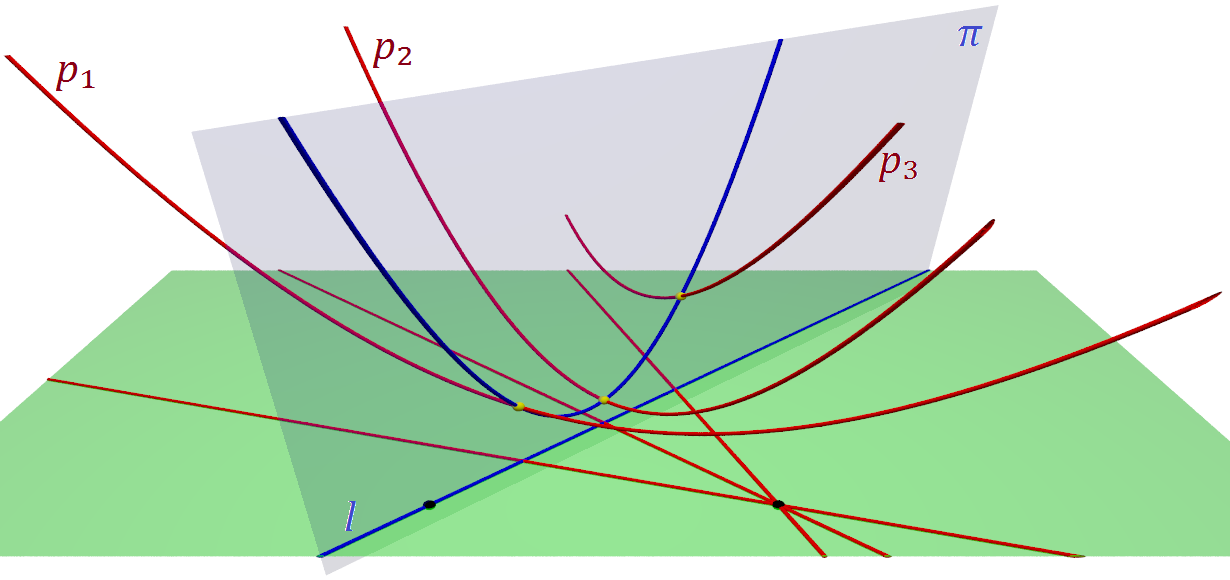}}
\end{minipage}
\caption{to Remark~\ref{rem:geom-constr-parab}: a geometric construction of surfaces containing two parabolas with vertical axes through each point. Three parabolas $p_1,p_2,p_3$ from the red family determine all parabolas from the blue family, i.~e., the whole surface.}
\label{fig:geom-constr-parab}
\end{figure}

\begin{remark}\label{rem:geom-constr-parab}
Surfaces described in Theorem~\ref{thm:main1} admit the following geometric construction (see Fig.~\ref{fig:geom-constr-parab} to the left). Let $\Omega$ be the set of all lines tangent to a fixed conic. Take three parabolas $p_1,p_2,p_3$ with the top views being lines lying in $\Omega$. For a generic line $l\in\Omega$ take the plane $\pi$ such that $\pi\perp Oxy$ and $\pi\cap Oxy=l$. Draw the parabola with a vertical axis through the points $p_1\cap\pi,p_2\cap\pi,p_3\cap\pi$. Then one can check by direct computation that the surface formed by all such parabolas for all $l\in\Omega$ contains two parabolas with vertical axes through each point.

The same construction still works if $\Omega$ is a union of two pencils of lines. In this case, the top views of $p_1,p_2,p_3$ should lie in one pencil and the line $l$ should be taken from the other one (see Fig.~\ref{fig:geom-constr-parab} to the right).
\end{remark}

\begin{theorem}\label{thm:main2}
Assume that through each point of an analytic surface in $\R^3$ one can draw two transversal arcs of isotropic circles fully contained in the surface (and analytically depending on the point). Assume that these two arcs lie neither in the same isotropic sphere nor in the same plane. Assume that through each point in some dense subset of the surface one can draw only finitely many (not nested) arcs of isotropic circles and line segments contained in the surface. Then the surface (possibly besides a one-dimensional subset) has a parametrization
\begin{equation}\label{eqn:mainparam2}
\Phi(u,v)=\left(\frac{P_0P_1-P_2P_3}{P_0^2+P_3^2},\frac{P_1P_3+P_0P_2}{P_0^2+P_3^2},\frac{Z}{P_0^2+P_3^2}\right)
\end{equation}
for some $P_0,P_1,P_2,P_3\in\R_{1,1}\subset\R[u,v]$ and $Z\in\R_{2,2}\subset\R[u,v]$, where $P_0^2+P_3^2\ne0$, such that the arcs of the isotropic circles are the curves $u=\const$ and $v=\const$.
\end{theorem}

\begin{figure}[h!]
\begin{minipage}{0.28\linewidth}
\center{\includegraphics[width=1\linewidth]{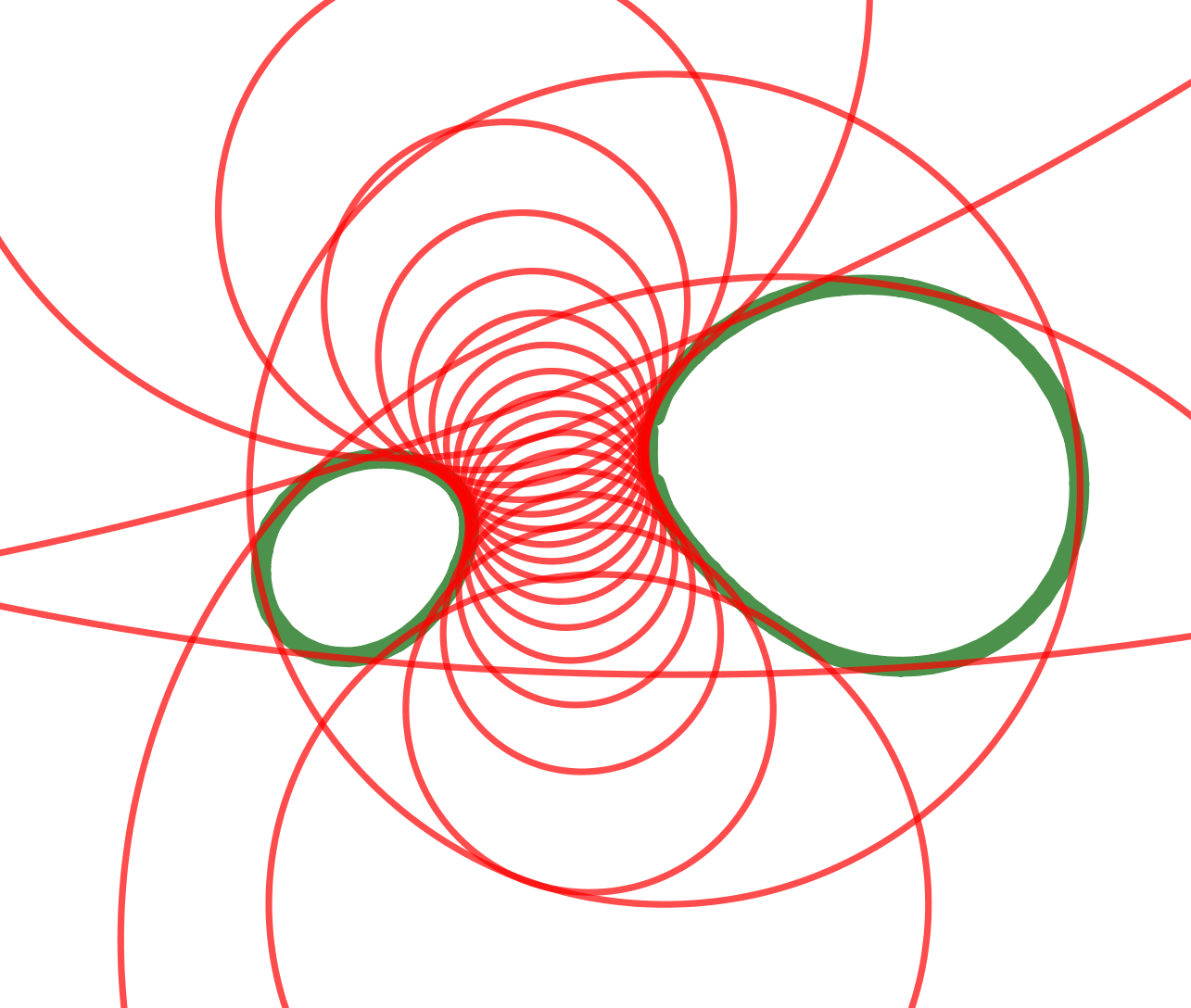}}
\end{minipage}
\hfill
\begin{minipage}{0.28\linewidth}
\center{\includegraphics[width=1\linewidth]{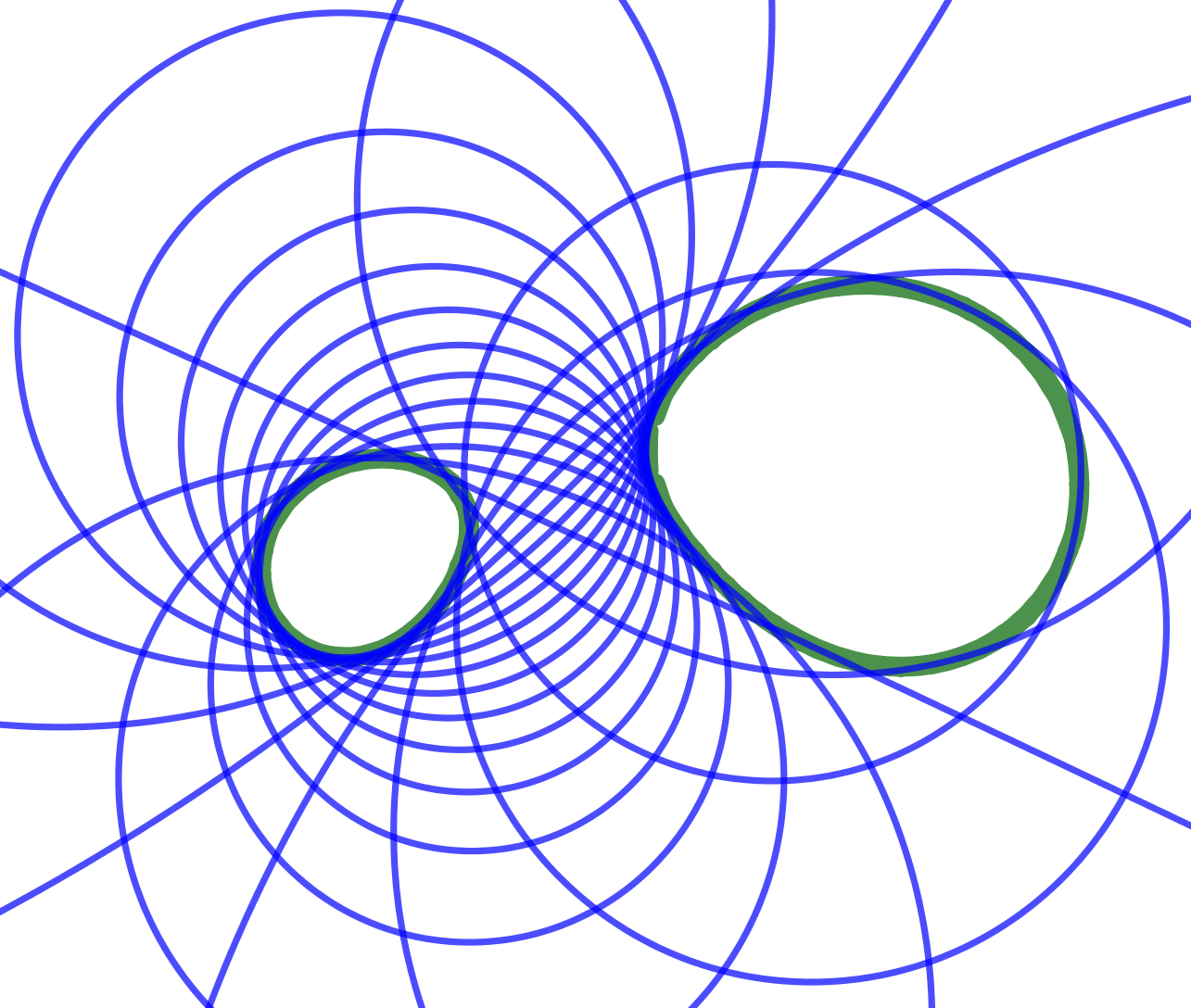}}
\end{minipage}
\vfill
\begin{minipage}{0.28\linewidth}
\center{\includegraphics[width=1\linewidth]{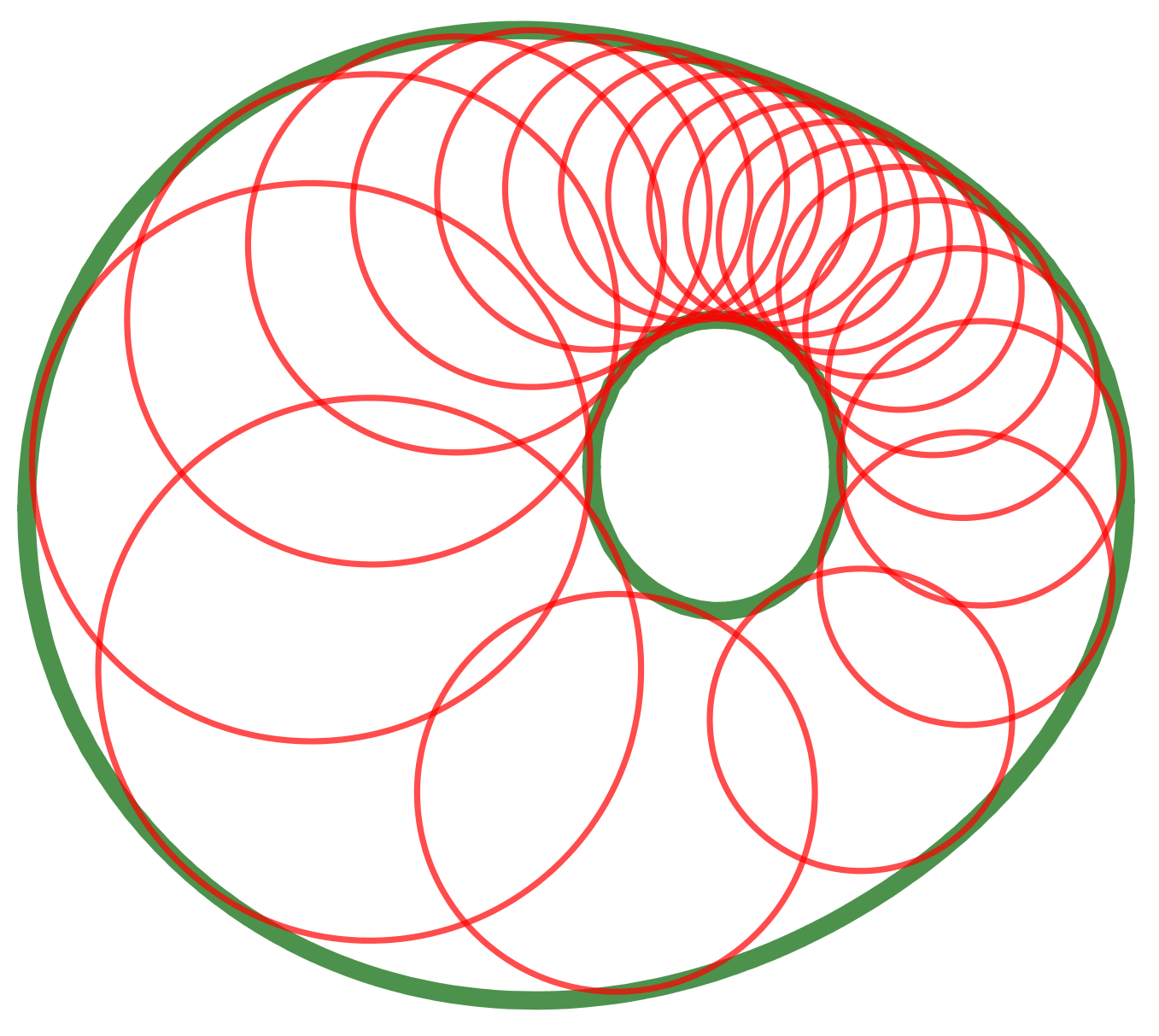}}
\end{minipage}
\hfill
\begin{minipage}{0.28\linewidth}
\center{\includegraphics[width=1\linewidth]{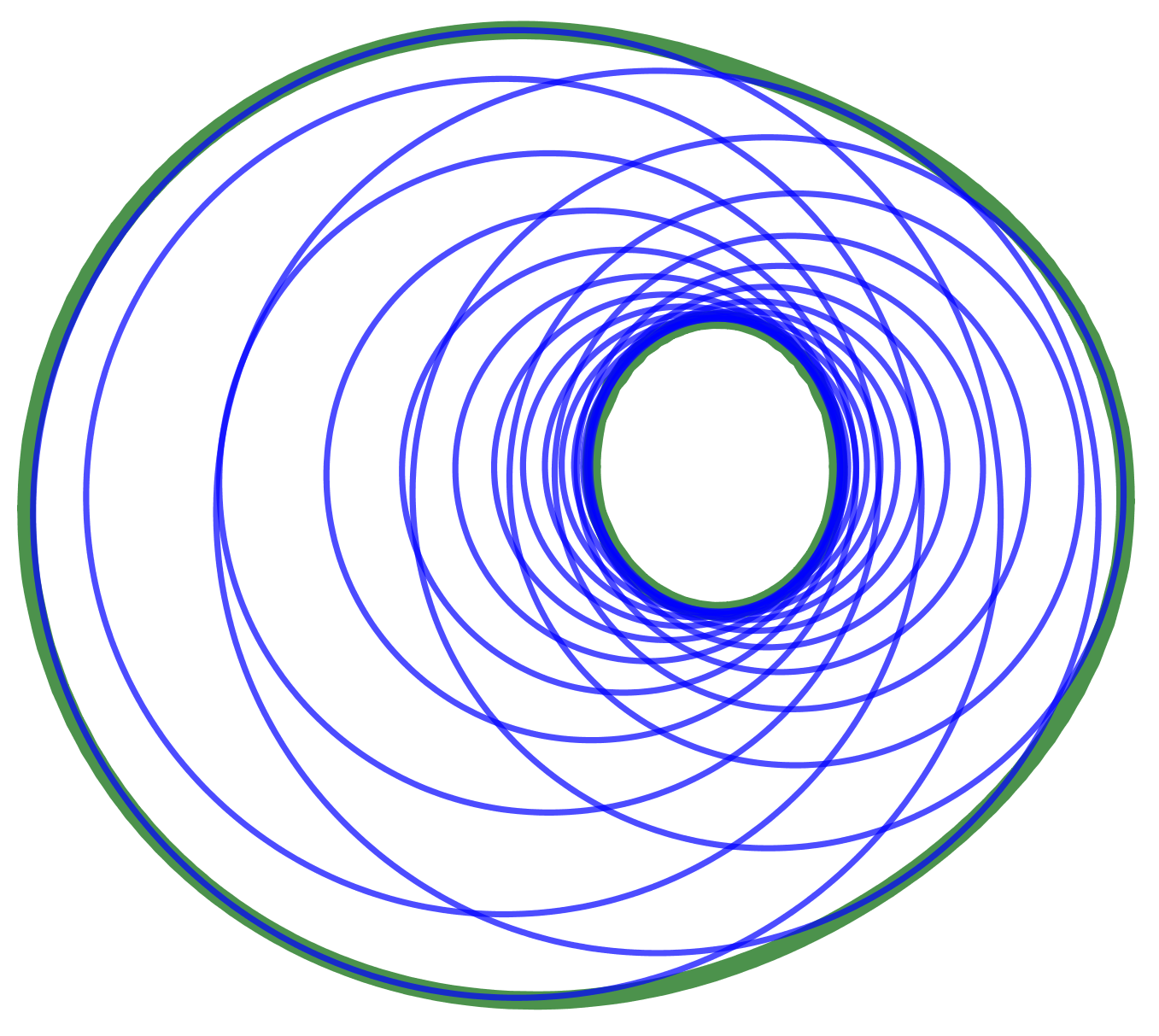}}
\end{minipage}
\caption{to Corollary~\ref{cor:top2}: top views of isotropic circles through each point. Top row: the top views of both families of isotropic circles have an envelope. Bottom row: the top view of the red family of isotropic circles has an envelope and the top view of the blue family has not. In both cases the envelope (green curve) is a cyclic.}
\label{fig:top2}
\end{figure}

\begin{corollary}\label{cor:top2}
Consider a surface satisfying the assumptions of Theorem~\ref{thm:main2}. Then if the top views of the isotropic circles $u=\const$ are all tangent to some regular smooth curve, then this curve is contained in a cyclic (see Fig.~\ref{fig:examples-2} and~\ref{fig:top2}). In this case, if the top views of the isotropic circles $v=\const$ are all tangent to some regular smooth curve as well, then this curve is contained in the same cyclic (see Fig.~\ref{fig:top2} top row).
\end{corollary}

\begin{remark}
Actually it is possible that one family of top views has a (real) envelope and the other one has not (see Fig.~\ref{fig:top2} bottom row). A more general assertion covering this case is stated in Conjecture~\ref{con:topview}.
\end{remark}

\begin{remark}
In fact, a surface given by a generic equation of form~\eqref{eqn:mainparam1} (respectively, \eqref{eqn:mainparam2}) contains two arcs of parabolas with vertical axes (respectively, two arcs of isotropic circles) through each point. 
\end{remark}

\begin{figure}[h]
\center{\includegraphics[width=0.5\linewidth]{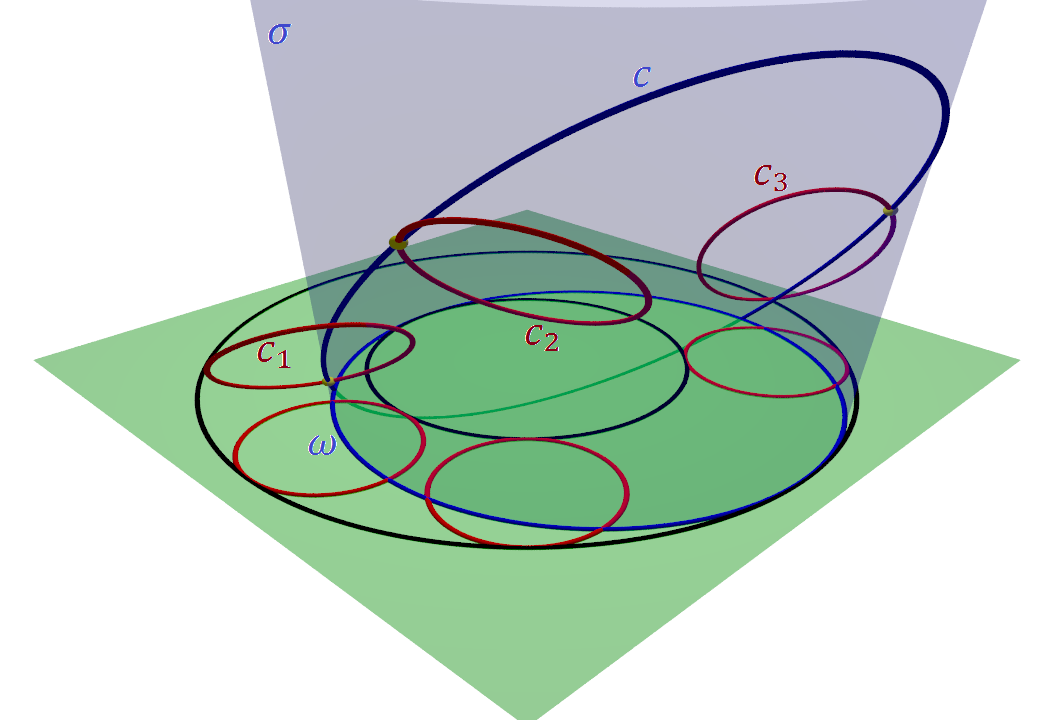}}
\caption{to Remark~\ref{rem:geom-constr-isot}: a geometric construction of surfaces containing two isotropic circles through each point. Three isotropic circles $c_1,c_2,c_3$ from the red family determine all isotropic circles from the blue family, i.~e., the whole surface.}
\label{fig:geom-constr-isot}
\end{figure}

\begin{remark}\label{rem:geom-constr-isot}
Surfaces described in Theorem~\ref{thm:main2} admit a geometric construction similar to the one described in Remark~\ref{rem:geom-constr-parab} (see Fig.~\ref{fig:geom-constr-isot}). Fix a cyclic together with two families of circles doubly tangent to the cyclic. Take three isotropic circles $c_1,c_2,c_3$ with top views lying in the first family. For a generic circle $\omega$ lying in the second family consider the right circular cylinder $\sigma$ such that $\sigma\cap Oxy=\omega$. Let $A_i$ be any of two points in the intersection $\sigma\cap c_i$ for $i=1,2,3$. The points $A_1,A_2,A_3$ determine a unique isotropic circle $c$. A continuous variation $\omega(t)$ of the circle $\omega$ in its family determines continuous variations of the points $A_1,A_2,A_3$ and so a continuous variation $c(t)$ of the isotropic circle $c$. We conjecture that the surface traced by the isotropic circles $c(t)$ for different $t$ contains two isotropic circular arcs through each point.
\end{remark}

\section{Outline of the proofs}\label{sec:outline}
We begin with some notation. Let $\mathbb{RP}^4$ be the 4-dimensional real projective space with the homogeneous coordinates ${(x_1:x_2:x_3:x_4:x_5)}$. Consider the following subsets:
\begin{itemize}
\item{the cylinder $S^2\times\R$ given by $x_1^2+x_2^2+x_4^2=x_5^2$ and $x_5\ne 0$;}
\item{the affine subspace $H_5$ given by $x_4=0,x_5=1$;}
\item{the affine line $l\subset S^2\times\R$ given by $x_1=x_2=0,x_4=x_5=1$.}
\end{itemize}

By a \emph{conic in $S^2\times\R$} we mean a nonempty section of $S^2\times\R$ by a 2-dimensional projective plane that contains no rulings of $S^2\times\R$ (so that any conic in $S^2\times\R$ is an irreducible plane conic).

The proofs of both Theorems~\ref{thm:main1} and~\ref{thm:main2} are similar and consist of three steps.
\begin{enumerate}
\item{Reducing the problem of finding all surfaces in $\R^3$ containing two isotropic circles through each point (or two parabolas through each point) to the problem of finding all surfaces in $S^2\times\R$ containing two conics through each point (in case of Theorem~\ref{thm:main1} these conics must intersect the line $l$).}
\item{Solving the resulting problems using parametrization of surfaces containing two conics through each point.}
\item{Extracting the solutions of the initial problems from the solutions obtained in the previous step.}
\end{enumerate}

For the first and the third steps we perform the inverse isotropic stereographic projection of a surface in $\R^3$ to obtain a surface in $S^2\times\R$. The \emph{isotropic stereographic projection} $\pi\colon (S^2\times\R)\setminus l\to H_5$ is defined by the formula
$$
\pi\colon (x_1:x_2:x_3:x_4:x_5)\mapsto\left(\frac{x_1}{x_5-x_4}:\frac{x_2}{x_5-x_4}:\frac{x_3}{x_5-x_4}:0:1\right).
$$
The inverse map $\pi^{-1}\colon H_5\to (S^2\times\R)\setminus l$ is given by
$$
\pi^{-1}\colon (x_1:x_2:x_3:0:1)\mapsto(2x_1:2x_2:2x_3:x_1^2+x_2^2-1:x_1^2+x_2^2+1).
$$

\begin{figure}[h]
\center{\includegraphics[scale=0.4]{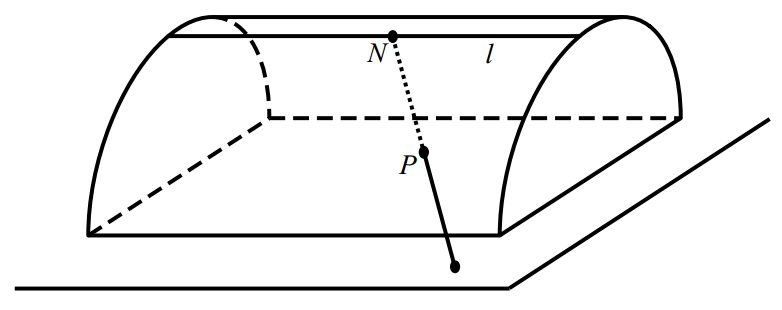}}
\caption{the isotropic stereographic projection (in 3-dimensional space).}
\label{fig:isot-stereo-proj}
\end{figure}

The point $N=(0:0:0:1:1)$ is called \emph{the projection center}. Thinking geometrically, the map $\pi$ takes a point $P\in(S^2\times\R)\setminus l$ to the intersection of the line $NP$ with $H_5$ (see Fig.~\ref{fig:isot-stereo-proj}). In particular, the image of any section of $(S^2\times\R)\setminus l$ by a projective subspace containing $N$ is an affine subspace of $H_5$.

The key property of the isotropic stereographic projection is stated in the following proposition, which is proved in~\cite{krasauskas2014bilinear}.
\begin{proposition}[{\cite[Theorem~3]{krasauskas2014bilinear}}]\label{prp:proj-prop}
Identify $H_5$ with $\R^3$ by the map $(x_1:x_2:x_3:0:1)\mapsto(x_1,x_2,x_3)$. Then the maps $\pi$ and $\pi^{-1}$ define a 1-1 correspondence between isotropic circles (respectively, isotropic spheres that are not cylinders) in $H_5$ and nonempty sections of $(S^2\times\R)\setminus l$ by a 2-dimensional (respectively, 3-dimensional) projective subspace that does not pass through the projection center $N$ and contains no rulings of $S^2\times\R$.
\end{proposition}

For the second step we need several auxiliary results. Our main instrument is the following theorem by J.\,Schicho. Recall that an \emph{analytic surface} in complex projective space ${P}^n=\mathbb{CP}^n$ is the image of an injective complex analytic map from a domain in $\C^2$ into ${P}^n$ with nondegenerate differential at each point. 

\begin{theorem}[{\cite[Theorem~1]{schicho2001multiple}}, {\cite[Theorem~4.1]{skopenkov2019surfaces}}]\label{thm:schicho}
Assume that through each point of an analytic surface $\Phi$ in a domain in ${P}^n$ one can draw two transversal conics intersecting each other only at this point (and analytically depending on the point) such that their intersections with the domain are contained in the surface. Assume that through each point in some dense subset of the surface one can draw only finitely many conics such that their intersections with the domain are contained in the surface. Then the surface is algebraic and has a parametrization (possibly besides a one-dimensional subset):
$$
\Phi(u,v)=(X_1(u,v):\dots:X_{n+1}(u,v)),
$$
for some $X_1,\dots,X_{n+1}\in\C_{2,2}$ such that the conics are the curves $u=\const$ and $v=\const$.
\end{theorem}

Using this theorem we find a convenient parametrization of surfaces in $S^2\times\R$ containing two conics through each point.

\begin{proposition}\label{prp:param}
Assume that through each point $P$ of an analytic surface $\Psi$ in $S^2\times\R$ one can draw two transversal arcs of conics fully contained in the surface (and analytically depending on the point). Assume that the conics containing these two arcs intersects only at $P$. Assume that through each point in some dense subset of the surface one can draw only finitely many (not nested) arcs of conics fully contained in the surface. Then the surface (possibly besides a one-dimensional subset) has a parametrization
\begin{equation}\label{eqn:cylinderparam}
\Psi(u,v)=(X_1(u,v):X_2(u,v):X_3(u,v):X_4(u,v):X_5(u,v)),
\end{equation}
for some $X_1,\dots,X_5\in\R_{2,2}$ satisfying the equation
\begin{equation}\label{eqn:cylinder}
X_1^2+X_2^2+X_4^2=X_5^2,
\end{equation}
and such that the conics are the curves $u=\const$ and $v=\const$.
\end{proposition}

The proof of this proposition is quite technical and goes along the lines of the proof of~\cite[Theorem~1.2]{skopenkov2019surfaces}. However the original proof of~\cite[Theorem~1.2]{skopenkov2019surfaces} had a minor gap. A correction is included as an appendix in the arxiv.org version~\cite{skopenkov2022surfaces} of that paper, where the proofs of both~\cite[Theorem~1.2]{skopenkov2019surfaces} and Proposition~\ref{prp:param} are presented (see~\cite[Remark~4.17]{skopenkov2022surfaces}). An official erratum is submitted.

The following proposition gives an algebraic interpretation for the condition that both isotropic circles through each point of the initial surface in $\R^3$ are parabolas. It is used in the proof of Theorem~\ref{thm:main1} only.

\begin{proposition}\label{prp:goodparam}
Suppose that a surface $\Psi$ in $S^2\times\R$ has parametrization~\eqref{eqn:cylinderparam}, where $X_1,\dots,X_5\in\R_{2,2}$ satisfy~\eqref{eqn:cylinder}. Assume that all the curves $u=\const$ and $v=\const$ on the surface $\Psi$ intersect the line $l$. Then the surface $\Psi$ has another parametrization (possibly besides a one-dimensional subset as well)
$$
\Psi(u,v)=(Y_1(u,v):\dots:Y_5(u,v)),
$$
for some $Y_1,\dots,Y_5\in\R_{2,2}$ satisfying~\eqref{eqn:cylinder} (with $X_i$ replaced by $Y_i$) such that $u=\const$ and $v=\const$ are the same curves and $Y_1,Y_2,Y_5-Y_4$ have a common divisor of degree at least~1 in $u$ and $v$.
\end{proposition}

The required parametrization~\eqref{eqn:mainparam1} in Theorem~\ref{thm:main1} is deduced from the following proposition.

\begin{proposition}\label{prp:polynomials}
Let $X_1,X_2,X_4,X_5\in\R_{2,2}$ be four polynomials satisfying~\eqref{eqn:cylinder}. Assume that $X_1$, $X_2$, and $X_5-X_4$ have a common divisor of degree at least~1 in $u$ and $v$. Then there exist polynomials $P,Q,R,T\in\R_{2,2}$ such that
\begin{align}
X_1&=2PRT,\notag\\
X_2&=2QRT,\label{eqn:tparam}\\
X_4&=(P^2+Q^2-R^2)T,\notag\\
X_5&=(P^2+Q^2+R^2)T.\notag
\end{align}
\end{proposition}

The required parametrization~\eqref{eqn:mainparam2} in Theorem~\ref{thm:main2} is deduced from the following proposition, which is proved in~\cite{dietz1993algebraic}.

\begin{proposition}[Parametrization of Pythagorean 4-tuples {\cite[Theorem~2.2]{dietz1993algebraic}}]\label{prp:4-tuples}
Let $X_1,X_2,X_4,X_5\in\R[u,v]$ be four polynomials satisfying~\eqref{eqn:cylinder}. Then there exist polynomials $P_0,P_1,P_2,P_3,T\in\R[u,v]$ such that 
\begin{align}
X_1&=2(P_0P_1-P_2P_3)T,\notag\\
X_2&=2(P_1P_3+P_0P_2)T,\label{eqn:4param}\\
X_4&=(P_1^2+P_2^2-P_0^2-P_3^2)T,\notag\\
X_5&=(P_0^2+P_1^2+P_2^2+P_3^2)T.\notag
\end{align}
\end{proposition}

\begin{remark}
Actually from~\cite[Theorem~2.2]{dietz1993algebraic} it follows that under the assumptions of Proposition~\ref{prp:4-tuples} there exist parametrization~\eqref{eqn:4param} with $X_5=\pm(P_0^2+P_1^2+P_2^2+P_3^2)T$. However one can exclude the case of minus sign in the latter formula by the change of variables $(P_0,P_1,P_2,P_3,T)\mapsto (-P_1,P_0,P_3,-P_2,-T)$.
\end{remark}

For the proof of Corollary~\ref{cor:top2} we need the following auxiliary results not pretending to be new. Consider a surface of form~\eqref{eqn:mainparam2}. Identify the plane $Oxy$ with $\C$ by the map $(x,y)\mapsto x+yi$. Then the projection of a point $\Phi(u,v)$ to the plane $Oxy$ is identified with the complex number
\begin{equation}\label{eqn:phi-bil-frac}
\frac{P_0P_1-P_2P_3}{P_0^2+P_3^2}+\frac{P_1P_3+P_0P_2}{P_0^2+P_3^2}\cdot i=\frac{P_1(u,v)+P_2(u,v)i}{P_0(u,v)-P_3(u,v)i}=:\frac{a_{11}uv+a_{10}u+a_{01}v+a_{00}}{b_{11}uv+b_{10}u+b_{01}v+b_{00}}
\end{equation}
for some $a_{ij},b_{ij}\in\C$. Therefore studying the top view of isotropic circles on $\Phi$ is closely related to classification of complex \emph{bilinear}-fractional maps of the form
\begin{equation}\label{eqn:bil-frac}
F(u,v)=\frac{a_{11}uv+a_{10}u+a_{01}v+a_{00}}{b_{11}uv+b_{10}u+b_{01}v+b_{00}}.
\end{equation}
This is the missing \emph{complex} analogue of the \emph{quaternionic} classification from~\cite[\S3]{skopenkov2019surfaces}, cf.~\cite{uhlig1976canonical}.

Two rational functions $F(u,v)$ and $G(u,v)$ of form~\eqref{eqn:bil-frac} are \emph{equivalent}, if there are invertible complex linear-fractional maps $f(z),f_u(u),f_v(v)$ such that
$$
F(u,v)=f(G(f_u(u),f_v(v))).
$$

\begin{theorem}[classification of complex bilinear-fractional maps]\label{thm:mp-reduction}
Each rational function $F(u,v)$ of form~\eqref{eqn:bil-frac} is equivalent to one of the following 5 polynomials: $uv,u+v,u,v,0$.
\end{theorem}

The latter theorem allows us to obtain the following description of the top view as a ``Minkowski sum or product'' of circles or lines. A \emph{generalized circle} in the plane is either a circle or a line. For a subset $A\subset\C$ and a point $z\in\C$ denote $z\cdot A=\{za\mid a\in A\}$ and $z+A=\{z+a\mid a\in A\}$. Note that if $A$ is a circle (respectively, a line) in the complex plane, then both sets $z\cdot A$ and $z+A$ are circles (respectively, lines) as well.

\begin{proposition}\label{prp:top2-ext}
Let $\Phi$ be a surface satisfying the assumptions of Theorem~\ref{thm:main2}. Then there exist generalized circles $\omega_1,\omega_2\subset\C$ and a linear-fractional map $f(z)$ such that
\begin{itemize}
\item[\emph{(i)}]{for any $u_0\in\R$ there exists $z_0\in\omega_1$ such that the top view of the arc of the isotropic circle $u=u_0$ on $\Phi$ is contained either in $\cl f(z_0\cdot\omega_2)$ or in $\cl f(z_0+\omega_2)$;}
\item[\emph{(ii)}]{for any $v_0\in\R$ there exists $w_0\in\omega_2$ such that the top view of the arc of the isotropic circle $v=v_0$ on $\Phi$ is contained either in $\cl f(w_0\cdot\omega_1)$ or in $\cl f(w_0+\omega_1)$.}
\end{itemize}
Here ``$\cl$'' stands for the closure of a set.
\end{proposition}

The last result required for the proof of Corollary~\ref{cor:top2} is the following one.

\begin{proposition}\label{prp:two-circles}
\begin{itemize}
\item[\emph{(i)}]{Let $\omega_1$ and $\omega_2$ be two generalized circles in the complex plane. Suppose that the generalized circles from the family $\Pi_1=\{w\cdot\omega_1\mid w\in\omega_2\}$ are tangent to some regular smooth curve; then this curve is contained in a cyclic. If, in addition, the generalized circles from the family $\Pi_2=\{z\cdot\omega_2\mid z\in\omega_1\}$ are tangent to some regular smooth curve as well, then the latter curve is contained in the same cyclic.}
\item[\emph{(ii)}]{The same assertions hold for the families $\Sigma_1=\{w+\omega_1\mid w\in\omega_2\}$ and $\Sigma_2=\{z+\omega_2\mid z\in\omega_1\}$.}
\end{itemize}
\end{proposition}

\section{Proofs}
\subsection{Proof of Proposition~\ref{prp:goodparam}.}\label{sec:goodparam}
\begin{lemma}\label{lem:roots}
Let $\Psi$ be a surface satisfying the assumptions of Proposition~\ref{prp:goodparam}. Then $\Psi$ has a parametrization (possibly besides a one-dimensional subset)
$$
\Psi(u,v)=(Y_1(u,v):\dots:Y_5(u,v)),
$$
for some $Y_1,\dots,Y_5\in\R_{2,2}$ satisfying~\eqref{eqn:cylinder} (with $X_i$ replaced by $Y_i$) such that $u=\const$ and $v=\const$ are the same curves as in Proposition~\ref{prp:goodparam} and
\begin{itemize}
\item[\emph{(i)}]{there are infinitely many $u_0\in\R$ such that the polynomials $Y_1(u_0,v)$, $Y_2(u_0,v)$, $Y_5(u_0,v)-Y_4(u_0,v)\in\R[v]$ have a common root $v_1$ (possibly depending on $u_0$);}
\item[\emph{(ii)}]{there are infinitely many $v_0\in\R$ such that the polynomials $Y_1(u,v_0)$, $Y_2(u,v_0)$, $Y_5(u,v_0)-Y_4(u,v_0)\in\R[u]$ have a common root $u_1$ (possibly depending on $v_0$).}
\end{itemize}
\end{lemma}

\noindent\textbf{Proof of Lemma~\ref{lem:roots}.} In parametrization~\eqref{eqn:cylinderparam}, the pair $(u,v)$ runs through some open subset of $\R^2$. Take any $(u_0,v_0)$ from this subset and consider the curve $u=u_0$ in $\Psi$. Since $\deg_v X_i\le 2$ for $i=1,\dots,5$, we obtain that this curve is either a smooth conic or a line. But there are no lines in $S^2\times\R$ intersecting the line $l$ except for the line $l$ itself. Hence the curve $u=u_0$ is either a smooth conic (in particular, $\max\limits_{i=1,\dots,5}\deg_v X_i=2$) or the line $l$. By the definition of $l$ it follows that two cases are possible:

\textit{Case 1: there exist $v_1,a\in\R$ such that}
$$
(X_1(u_0,v_1):\dots:X_5(u_0,v_1))=(0:0:a:1:1).
$$

\textit{Case 2: for some $a\in\R$ we have}
$$
(v^2X_1(u_0,1/v):\dots:v^2X_5(u_0,1/v))|_{v=0}=(0:0:a:1:1).
$$

If for infinitely many $u_0$ we have the first case, then condition~(i) of Lemma~\ref{lem:roots} holds (if we take $Y_i=X_i$ for $i=1,\dots,5$). Otherwise for infinitely many $u_0$ we have the second case. Then let us change the parametrization taking $X_i'(u,v)=v^2X_i(u,1/v)$ for $i=1,\dots,5$. For the resulting parametrization $(X_1'(u,v):\dots:X_5'(u,v))$ condition~(i) holds. Repeating this procedure with the curve $v=v_0$ we obtain a parametrization $(Y_1(u,v):\dots:Y_5(u,v))$ which satisfies both conditions~(i) and~(ii).\qed

\begin{lemma}\label{lem:gcd}
Let $F_1$, $F_2$, and $F_3$ be three polynomials in $\R[u,v]$. Suppose that for infinitely many $u_0\in\R$ the polynomials $F_1(u_0,v)$, $F_2(u_0,v)$, and $F_3(u_0,v)$ have a common root. Then $F_1$, $F_2$, and $F_3$ have a common divisor of degree at least one in $v$.
\end{lemma}

\noindent\textbf{Proof of Lemma~\ref{lem:gcd}.} Let $D$ be one of the greatest common divisors of $F_1$, $F_2$, and $F_3$. First let us prove that $\deg D\ge1$. Without loss of generality it can be assumed that $F_1\ne0$. For the case $F_2=F_3=0$ there is nothing to prove. Otherwise without loss of generality it can be assumed that $F_2\ne0$. Then the curves $F_1(u,v)=0$ and $F_2(u,v)=0$ have infinitely many common points. By Bezout's theorem it follows that $F_1$ and $F_2$ are not coprime. Let $A$ be one of the greatest common divisors of $F_1$ and $F_2$, so that $F_1=AF_1'$, $F_2=AF_2'$ for some $F_1',F_2'\in\R[u,v]$ and $\deg A\ge1$. If $F_3=0$, then $A$ is a common divisor of $F_1$, $F_2$, and $F_3$; in particular, $A\mid D$ and we are done. Otherwise the curves $F_1=0$, $F_2=0$, and $F_3=0$ have infinitely many common points, whereas the curves $F_1'=0$ and $F_2'=0$ have finitely many common points. Hence the curves $A=0$ and $F_3=0$ have infinitely many common points. By Bezout's theorem it follows that $A$ and $F_3$ are not coprime. Thus $F_1$, $F_2$, and $F_3$ are not coprime as well.

Further, let us prove that actually $\deg_v D\ge1$. Assume the converse. Then $D\in\R[u]$. Suppose $F_i=DF_i''$ for some $F_i''\in\R[u,v]$ (here $i=1,2,3$). Then there exist infinitely many $u_0$ such that $D(u_0)\ne0$ but the polynomials $F_1(u_0,v)$, $F_2(u_0,v)$, and $F_3(u_0,v)$ have a common root. For any such $u_0$ the polynomials $F_1''(u_0,v)$, $F_2''(u_0,v)$, and $F_3''(u_0,v)$ have a common root. From the previous paragraph it follows that $F_1''$, $F_2''$, and $F_3''$ are not coprime. This contradicts the assumption that $D$ is the greatest common divisor of $F_1$, $F_2$, and $F_3$, which completes the proof.\qed

\noindent\textbf{Proof of Proposition~\ref{prp:goodparam}} follows from Lemmas~\ref{lem:roots} and~\ref{lem:gcd}.\qed

\subsection{Proof of Proposition~\ref{prp:polynomials}.}\label{sec:polynomials}
First let us prove the following simple lemma.

\begin{lemma}\label{lem:sos}
Let $F_1$ and $F_2$ be two coprime polynomials in $\R[u,v]$ (in particular, $F_1^2+F_2^2\ne0$). Suppose that a real polynomial $I$ divides $F_1^2+F_2^2$ and $I\ne\const$; then $\deg_u I\ge2$ or $\deg_v I\ge2$.
\end{lemma}
\noindent\textbf{Proof of Lemma~\ref{lem:sos}.} Assume the converse. Without loss of generality it can be assumed that $\deg_v I=1$. Then there exist infinitely many $u_0\in\R$ such that $I(u_0,v)$ has a real root $v_0$. For any such pair $(u_0,v_0)$ we have $F_1^2(u_0,v_0)+F_2^2(u_0,v_0)=0$, i.~e., $F_1(u_0,v_0)=F_2(u_0,v_0)=0$. In particular, the curves $F_1=0$ and $F_2=0$ have infinitely many common points. From Bezout's theorem it follows that $F_1$ and $F_2$ are not coprime. This contradiction concludes the proof.\qed

\smallskip

\noindent\textbf{Proof of Proposition~\ref{prp:polynomials}.} If $X_1$ and $X_2$ vanish simultaneously, then by \eqref{eqn:cylinder} we have $X_4=\pm X_5$ and the required parametrization is given by either $P=Q=0,R=1,T=X_5$ (in the case when $X_4=-X_5$) or $Q=R=0,P\hm=1,T=X_4$ (in the case when $X_4=X_5$). Assume further that $X_1^2+X_2^2\ne0$.

Let $D$ be one of the greatest common divisors of $X_1$, $X_2$, and $X_5-X_4$. Denote $X_1=DP$, $X_2=DQ$, $X_5-X_4=DR$. By~\eqref{eqn:cylinder} we have
\begin{equation}\label{eqn:dpqrx}
D(P^2+Q^2)=R(X_5+X_4).
\end{equation}

Let us show that $R\mid D$; then taking $T=D/2R$ we obtain the required parametrization. By~\eqref{eqn:dpqrx} it suffices to show that $R$ and $P^2+Q^2$ are coprime. Since $X_1^2+X_2^2\ne0$, we have $P^2+Q^2\ne0$. Hence $P$ and $Q$ have some greatest common divisor $D'$ so that $P=D'P'$, $Q=D'Q'$, where $P'$ and $Q'$ are coprime. Then $R$ and $D'$ are coprime because otherwise $D$ is not a greatest common divisor of $X_1$, $X_2$, and $X_5-X_4$.

It remains to show that $R$ and $(P')^2+(Q')^2$ are coprime. By the assumptions of Proposition~\ref{prp:polynomials} we have $\deg_v D\ge1$ and $\deg_u D\ge1$. However $DR=X_5-X_4\in\R_{2,2}$, hence $R\in\R_{1,1}$. Let $I$ be an arbitrary irreducible (in $\R[u,v]$) nonconstant divisor of $(P')^2+(Q')^2$. By Lemma~\ref{lem:sos} it follows that $\deg_u I\ge2$ or $\deg_v I\ge2$. Since $R\in\R_{1,1}$, it follows that $I\nmid R$. Thus $R$ and $(P')^2+(Q')^2$ are coprime and $R\mid D$.\qed

\subsection{Proof of Theorem~\ref{thm:main1}.}\label{sec:main1}
Identify space $\R^3$ with $H_5$ by the map $(x,y,z)\mapsto(x:y:z:0:1)$. Consider the image $\Psi$ of the given surface $\Phi$ under the inverse isotropic stereographic projection $\pi^{-1}$.

Let us show that $\Psi$ satisfies the assumptions of Proposition~\ref{prp:param}. By Proposition~\ref{prp:proj-prop} it follows that each isotropic circle in $H_5$ maps to a conic (in particular, by the assumptions of Theorem~\ref{thm:main1} it follows that through each point of some dense subset of $\Psi$ one can draw only finitely many not nested arcs of conics). Moreover, each parabola with vertical axis is mapped to a conic that intersects the line $l$ (with the point on $l$ excluded). Indeed, the projective closure of a parabola with a vertical axis passes through the point $A=(0:0:1:0:0)$ while the line joining this point with the projection center $N$ is $l\cup\{A\}$. It follows that the closure of the image of a parabola intersects $l\cup\{A\}$ and the intersection point is not $A$ because otherwise the image of the parabola contains a ruling of $S^2\times\R$, which contradicts Proposition~\ref{prp:proj-prop}.

It remains to check the condition that the two conics containing the two arcs of conics through a point $P$ of $\Psi$ intersect only at $P$. Indeed, otherwise both conics are contained in one 3-dimensional subspace. If this subspace does not contain $N$ and rulings of $\S^2\times\R$, then by Proposition~\ref{prp:proj-prop} it follows that the $\pi$-images of the conics are contained in one isotropic sphere, which is forbidden. If this subspace contains $N$ or a ruling of $\S^2\times\R$, then the $\pi$-images of the conics are contained either in a plane or in a circular cylinder with vertical axis, which is forbidden as well.

Then by Propositions~\ref{prp:param},~\ref{prp:goodparam}, and~\ref{prp:polynomials} the surface $\Psi$ has parametrization~\eqref{eqn:cylinderparam}, where $X_1,\dots,X_5$ satisfy~\eqref{eqn:tparam}. We have $R\ne0$, otherwise $X_1=X_2=0$, $X_4=X_5$, and the surface is contained in a line, which is impossible.

Now, using the stereographic projection $\pi$, we obtain the parametrization of the initial surface $\Phi$:
$$
\Phi(u,v)=\left(\frac{P}{R},\frac{Q}{R},\frac{Z}{R^2T}\right)
$$
with $Z:=X_3/2$, for which the parabolas through each point are the curves $u=\const$ and $v=\const$.

Let us show that actually $T=\const$. Assume the converse. Then without loss of generality we may assume that $\deg_u T\ge1$. Then (in the notation of Proposition~\ref{prp:polynomials}) from $X_5=(P^2+Q^2+R^2)T$ and the condition $X_5\in\R_{2,2}$ we obtain $P,Q,R\in\R_{0,1}$. This means that all points of $\Phi$ with $v=\const$ are contained in a vertical line, i.~e., the curve $v=\const$ is not a parabola with vertical axis. This contradiction implies that $T=\const$. Since $X_5\in\R_{2,2}$ and $X_5=(P^2+Q^2+R^2)T$, we obtain $P,Q,R\in\R_{1,1}$. Absorbing the constant $T$ into $Z(u,v)$ we arrive at parametrization~\eqref{eqn:mainparam1}.\qed


\subsection{Proof of Theorem~\ref{thm:main2}}\label{sec:main2}
As in the proof of Theorem~\ref{thm:main1}, identify space $\R^3$ with $H_5$ by the map $(x,y,z)\mapsto(x:y:z:0:1)$ and consider the image $\Psi$ of the given surface $\Phi$ under inverse isotropic stereographic projection $\pi^{-1}$.

Let us show that $\Psi$ satisfies the assumptions of Proposition~\ref{prp:param}. By Proposition~\ref{prp:proj-prop} it follows that each isotropic circle in $H_5$ is mapped to a conic under this projection (in particular, by the assumptions of Theorem~\ref{thm:main2} it follows that through each point of some dense subset of $\Psi$ one can draw only finitely many not nested arcs of conics). Literally as in the 3rd paragraph of the proof of Theorem~\ref{thm:main1} it is shown that the two conics containing the two arcs through a point $P$ of $\Psi$ intersect only at $P$.

By Propositions~\ref{prp:param} and~\ref{prp:4-tuples} the surface $\Psi$ has parametrization~\eqref{eqn:cylinderparam}, where $X_1,\dots,X_5$ now satisfy~\eqref{eqn:4param}. We have $P_0^2+P_3^2\ne0$ and $P_1^2+P_2^2\ne0$, otherwise $X_1=X_2=0$, $X_4=\pm X_5$, and the surface is contained in a line, which is impossible.

Using the stereographic projection $\pi$, we obtain the parametrization of the initial surface $\Phi$:
$$
\Phi(u,v)=\left(\frac{P_0P_1-P_2P_3}{P_0^2+P_3^2},\frac{P_1P_3+P_0P_2}{P_0^2+P_3^2},\frac{Z}{(P_0^2+P_3^2)T}\right)
$$
with $Z:=X_3/2$, for which the isotropic circles through each point are the curves $u=\const$ and $v=\const$.

Let us show that $T=\const$. Assume the converse. Then without loss of generality we may assume that $\deg_u T\ge1$. Then (in the notation of Proposition~\ref{prp:4-tuples}) from $X_5=(P_0^2+P_1^2+P_2^2+P_3^2)T$ and the condition $X_5\in\R_{2,2}$ we obtain $P_i\in\R_{0,1}$ for $i=1,2,3,4$. This means that all points of $\Phi$ with $v=\const$ are contained in a vertical line, i.~e., the curve $v=\const$ is not an isotropic circle. This contradiction implies that $T=\const$. Since $X_5\in\R_{2,2}$ and $X_5=(P_0^2+P_1^2+P_2^2+P_3^2)T$, we obtain $P_i\in\R_{1,1}$ for $i=1,2,3,4$. Absorbing the constant $T$ into $Z(u,v)$ we arrive at parametrization~\eqref{eqn:mainparam2}.\qed

\subsection{Proof of Corollary~\ref{cor:top1}}\label{sec:top1}
%

Let us start from parametrization~\eqref{eqn:mainparam1} given by Theorem~\ref{thm:main1}. Identify our space $\R^3$ with $\{(x:y:z:t)\mid t\ne0\}\subset\mathbb{RP}^3$ by the map $(x,y,z)\mapsto(x:y:z:1)$. Write
$$
P(u,v)=p_{11}uv+p_{10}u+p_{01}v+p_{00},\;
Q(u,v)=q_{11}uv+q_{10}u+q_{01}v+q_{00},\;
R(u,v)=r_{11}uv+r_{10}u+r_{01}v+r_{00}.
$$

The top views of the parabolas through a point $\Phi(u_0,v_0)$ of the surface $\Phi$ are the two lines parametrized by
$$
(P(u,v_0):Q(u,v_0):0:R(u,v_0))\mbox{ and }(P(u_0,v):Q(u_0,v):0:R(u_0,v)).
$$

These lines are the images of the rulings $u/w=u_0,v/w=v_0$ of the quadric $\Sigma=\{(s:u:v:w)\mid sw=uv\}$ under the projective map
$$
(s:u:v:w)\mapsto(p_{11}s+p_{10}u+p_{01}v+p_{00}w:q_{11}s+q_{10}u+q_{01}v+q_{00}w:r_{11}s+r_{10}u+r_{01}v+r_{00}w:0).
$$
Now the corollary follows from the following well-known assertion: \emph{under a projective map $\mathbb{RP}^3\dashrightarrow\mathbb{RP}^2$ the rulings of a ruled quadric are mapped to either the tangents to a conic, or two pencils of lines, or one line, or one point.}\qed

\subsection{Proof of Theorem~\ref{thm:mp-reduction}.}\label{sec:linal}
Let $\mattc$ be the set of all ($2\times2$)-matrices with entries in $\C$. Take $M=\bigl(\begin{smallmatrix}a & b\\c & d\end{smallmatrix}\bigr)\in\gltc,A=\bigl(\begin{smallmatrix}a_{11} & a_{10}\\a_{01} & a_{00}\end{smallmatrix}\bigr)\in\mattc,B=\bigl(\begin{smallmatrix}b_{11} & b_{10}\\b_{01} & b_{00}\end{smallmatrix}\bigr)\in\mattc\setminus\{0\}$. By definition, put

\begin{gather*}
f_M(z):=\frac{az+b}{cz+d},\\
F^A_B(u,v):=\frac{a_{11}uv+a_{10}u+a_{01}v+a_{00}}{b_{11}uv+b_{10}u+b_{01}v+b_{00}}.
\end{gather*}

\begin{lemma}\label{lem:actions}
For any $M,N\in\gltc$ we have $f_N\circ f_M=f_{NM}$. For any $A\in\mattc,B\in\mattc\setminus\{0\}$ and $C,D\in\gltc$ we have $F^A_B(f_C(u),f_D(v))=F^{C^TAD}_{C^TBD}(u,v)$.
\end{lemma}

\begin{lemma}\label{lem:multable}
Let $P(u,v)=c_{11}uv+c_{10}u+c_{01}v+c_{00}\in\R_{1,1}$; then $P=QR$ for some $Q\in\R_{1,0},R\in\R_{0,1}$ if and only if $c_{00}c_{11}-c_{10}c_{01}=0$.
\end{lemma}

\noindent Proofs of Lemmas~\ref{lem:actions} and~\ref{lem:multable} can be obtained by a direct computation.

\smallskip

\noindent\textbf{Proof of Theorem~\ref{thm:mp-reduction}.} Suppose that $F(u,v)=F^A_B(u,v)$ for some $A\in\mattc,B\in\mattc\setminus\{0\}$.

\textit{Case 1: $\det A=\det B=0$.} Then by Lemma~\ref{lem:multable} we have
$$
F^A_B(u,v)=\frac{a_{11}uv+a_{10}u+a_{01}v+a_{00}}{b_{11}uv+b_{10}u+b_{01}v+b_{00}}=\frac{au+b}{cu+d}\cdot\frac{a'v+b'}{c'v+d'}=:f_u(u)\cdot f_v(v)
$$
for some $a,b,c,d,a',b',c',d'\in\C$. If $f_u(u),f_v(v)\ne\const$, then this shows that $F^A_B(u,v)$ is equivalent to $uv$. Otherwise $F^A_B(u,v)$ is clearly equivalent to either $u$, or $v$, or $0$.

\textit{Case 2: $\det B\ne0$.} Let us show that there exist matrices $C,D,M\in\gltc$ such that $f_M(F^A_B(f_C(u),f_D(v)))$ equals one of the polynomials $uv,u+v,u,v,0$. From Lemma~\ref{lem:actions} it follows that $f_M(F^A_B(f_C(u),f_D(v)))=f_M(F^{C^TAD}_{C^TBD}(u,v))$. Suppose that $X\in\gltc$ is such that $J=XAB^{-1}X^{-1}$ is the Jordan normal form of the matrix $AB^{-1}$. Consider the following 3 subcases, depending on the Jordan cells of the matrix $J$.

\textit{Subcase 2.1: $J=\bigl(\begin{smallmatrix}\lambda & 0\\0 & \mu\end{smallmatrix}\bigr),\lambda\ne\mu$.} Then take $C=X^T,D=B^{-1}X^{-1},M=\bigl(\begin{smallmatrix}-1 & \mu\\1 & -\lambda\end{smallmatrix}\bigr)$. We have
$$
f_M(F^{C^TAD}_{C^TBD}(u,v))=f_M(F^J_{\id}(u,v))=f_M\left(\frac{\lambda uv+\mu}{uv+1}\right)=uv.
$$

\textit{Subcase 2.2: $J=\bigl(\begin{smallmatrix}\lambda & 0\\0 & \lambda\end{smallmatrix}\bigr)$ (that is, $F^A_B(u,v)=\const$).} Then take $C=X^T,D=B^{-1}X^{-1},M=\bigl(\begin{smallmatrix}1 & -\lambda\\0 & 1\end{smallmatrix}\bigr)$. We have
$$
f_M(F^{C^TAD}_{C^TBD}(u,v))=f_M(F^J_{\id}(u,v))=f_M\left(\frac{\lambda uv+\lambda}{uv+1}\right)=0.
$$

\textit{Subcase 2.3: $J=\bigl(\begin{smallmatrix}\lambda & 1\\0 & \lambda\end{smallmatrix}\bigr)$.} Then take $C=X^TY,D=B^{-1}X^{-1},M=\bigl(\begin{smallmatrix}0 & 1\\1 & -\lambda\end{smallmatrix}\bigr)$, where $Y=\bigl(\begin{smallmatrix}0 & 1\\1 & 0\end{smallmatrix}\bigr)$. We have
$$
f_M(F^{C^TAD}_{C^TBD}(u,v))=f_M(F^{YJ}_Y(u,v))=f_M\left(\frac{\lambda(u+v)+1}{u+v}\right)=u+v.
$$

\textit{Case 3: $\det A\ne0,\det B=0$.} Note that $F^A_B(u,v)$ is equivalent to $F^B_A(u,v)$ because $F^B_A(u,v)=f_Y(F^A_B(u,v))$. Thus this case reduces to the previous one.\qed

\subsection{Proof of Proposition~\ref{prp:top2-ext}}
Let us start from parametrization~\eqref{eqn:mainparam2} given by Theorem~\ref{thm:main2}. The projection of a point $\Phi(u,v)$ to the plane $Oxy$ is given by formula~\eqref{eqn:phi-bil-frac}. The top views of the isotropic circles through a point $\Phi(u_0,v_0)$ are contained in the generalized circles $\cl\,\{F(u_0,v)\mid v\in\R\}$ and $\cl\,\{F(u,v_0)\mid u\in\R\}$, where $F(u,v)$ is given by~\eqref{eqn:bil-frac}.

By Theorem~\ref{thm:mp-reduction} there exist complex linear-fractional maps $f,f_z,f_w$ such that $G(z,w):=f(F(f_z(z),f_w(w)))$ equals one of the polynomials $zw,z+w,z,w,0$. Consider the following three cases.

\textit{Case 1: $G(z,w)=zw$.} Denote $\omega_1=\cl f^{-1}_z(\R)$, $\omega_2=\cl f^{-1}_w(\R)$, and also $z_0=f^{-1}_z(u_0),w_0=f^{-1}_w(v_0)$. Then
$$
\cl\,\{F(u_0,v)\mid v\in\R\}=\cl\,\{f^{-1}(f(F(f_z(z_0),f_w(w))))\mid w\in\omega_2\}=\cl\,\{f^{-1}(G(z_0,w))\mid w\in\omega_2\}=\cl f^{-1}(z_0\cdot\omega_2).
$$
Analogously, $\cl\,\{F(u,v_0)\mid u\in\R\}=\cl f^{-1}(w_0\cdot\omega_1)$.

\textit{Case 2: $G(z,w)=z+w$} is completely analogous to the previous case, only the product is replaced by the sum.

\textit{Case 3: $G(z,w)$ equals either $z$, or $w$, or $0$.} Then $G(z,w)$, hence $F(u,v)$, does not depend on one of the variables. Thus the top view of one of the isotropic circles $u=\const$ or $v=\const$ is a single point, which is impossible.\qed

\subsection{Proofs of Proposition~\ref{prp:two-circles} and Corollary~\ref{cor:top2}}\label{sec:top2}
Recall that any generalized circle in the complex plane can be given by an equation
\begin{equation}\label{eqn:circle}
\alpha z\bar z+\beta z+\bar\beta\bar z+\gamma=0,
\end{equation}
where $\alpha,\gamma\in\R,\beta\in\C$. Suppose that $P(z,\bar z)$ and $Q(z,\bar z)$ are two nonproportional equations of form~\eqref{eqn:circle}; then the family of generalized circles $\{\lambda P(z,\bar z)+\mu Q(z,\bar z)=0\mid\lambda,\mu\in\R\}$ is \emph{a pencil of circles}. Equivalently, a pencil of circles can be defined as the image of either a pencil of lines or the set of circles $\{|z|=R\mid R\in\R\}$ under a linear-fractional map. In particular, a pencil of circles never has a smooth envelope.

\smallskip

\noindent\textbf{Proof of Proposition~\ref{prp:two-circles}. (i)} (A.\,A.\,Zaslavsky, private communication). Let us show that the envelope of the family $\Pi_1$ is contained in a cyclic. Suppose that $\omega_2=\{\frac{av+b}{cv+d}\mid v\in\R\setminus\{-\frac{d}{c}\}\}\cup\{\frac{a}{c}\}$, where $a,b,c,d\in\C$ are fixed, and that $\omega_1$ is given by~\eqref{eqn:circle}. Then the generalized circle $\frac{av+b}{cv+d}\cdot\omega_1\in\Pi_1$ has equation
\begin{equation}\label{eqn:circles-long}
\alpha z\bar z(cv+d)(\bar cv+\bar d)+
\beta z(\bar av+\bar b)(cv+d)+
\bar\beta\bar z(av+b)(\bar cv+\bar d)+
\gamma(av+b)(\bar av+\bar b)=0.
\end{equation}
It can be rewritten as
\begin{equation}\label{eqn:circles-short}
A(z,\bar z)v^2+B(z,\bar z)v+C(z,\bar z)=0,
\end{equation}
where
\begin{align*}
A(z,\bar z)&=\alpha c\bar cz\bar z+\beta\bar acz+\bar\beta a\bar c\bar z+\gamma a\bar a,\\
B(z,\bar z)&=\alpha(c\bar d+\bar cd)z\bar z+\beta(\bar bc+\bar ad)z+\bar\beta(b\bar c+a\bar d)\bar z+\gamma(\bar ab+a\bar b),\\
C(z,\bar z)&=\alpha d\bar dz\bar z+\beta\bar bdz+\bar\beta b\bar d\bar z+\gamma b\bar b.
\end{align*}
In other words, $\Pi_1$ is a \emph{quadratic} family of generalized circles. To find the equation of the envelope, differentiate~\eqref{eqn:circles-short} with respect to $v$ and get
\begin{equation}\label{eqn:circles-dv}
2A(z,\bar z)v+B(z,\bar z)=0.
\end{equation}
The system of equations~\eqref{eqn:circles-short} and~\eqref{eqn:circles-dv} gives the envelope. Eliminating $v$ from the system, we get
\begin{equation}\label{eqn:the-cyclic}
B^2(z,\bar z)-4A(z,\bar z)C(z,\bar z)=0.
\end{equation}
It is easy to see that the substitution of $z=x+yi,\bar z=x-yi$ into~\eqref{eqn:the-cyclic} gives an equation in $x$ and $y$ of form~\eqref{eqn:cyclic-general}, i.~e., a cyclic.

Let us check that equation~\eqref{eqn:the-cyclic} does not hold identically. Indeed, otherwise we have $(2Av+B)^2=4A(Av^2+Bv+C)$, i.~e., for each $v\in\R$ the generalized circle given by~\eqref{eqn:circles-short} lies in the pencil of circles given by $2A\lambda+B\mu=0$ for $\lambda,\mu\in\R$. Hence the family $\Pi_1$ has no envelope, which contradicts to the assumptions of the proposition.

Now suppose that $\omega_1$ and $\omega_2$ are parametrized by some parameters $s_1\in\R$ and $s_2\in\R$ respectively. Then the family $\Pi_1$ can be considered as parametrized by $s_2$. Moreover, each curve of this family is itself parametrized by $s_1$. Changing the order of parameters, we obtain the family $\Pi_2$. If both families have a smooth envelope, then these envelopes coincide due to a general fact, see~\cite{green1952envelope}.

\noindent\textbf{(ii)} In this case the envelope of the family $\Sigma_1$ is either a pair of concentric circles, or a circle, or a pair of parallel lines (the proof is trivial). All these sets are cyclics.\qed

\smallskip

\noindent\textbf{Proof of Corollary~\ref{cor:top2}} follows from Propositions~\ref{prp:top2-ext} and~\ref{prp:two-circles} and the fact that cyclics are mapped to cyclics under complex linear-fractional transformations.\qed

\section{Open problems}\label{sec:open}
We believe that the finiteness of number of conics through each point in our main theorems can be dropped.
\begin{conjecture}
Assume that through each point of an analytic surface in $\R^3$ one can draw two transversal parabolic arcs with vertical axes fully contained in the surface (and analytically depending on the point). Assume that these two arcs lie neither in the same isotropic sphere nor in the same plane. Then the surface (possibly besides a one-dimensional subset) has a parametrization~\eqref{eqn:mainparam1}.
\end{conjecture}

\begin{conjecture}\label{con:main2}
Assume that through each point of an analytic surface in $\R^3$ one can draw two transversal arcs of isotropic circles fully contained in the surface (and analytically depending on the point). Assume that these two arcs lie neither in the same isotropic sphere nor in the same plane. Then the surface (possibly besides a one-dimensional subset) has a parametrization~\eqref{eqn:mainparam2}.
\end{conjecture}

The following conjecture is the isotropic analogue of~\cite[Theorem~3.5]{nilov2013surface}.
\begin{conjecture}
Assume that through each point of an analytic surface in $\R^3$ one can draw two transversal arcs of isotropic circles lying on one isotropic sphere and fully contained in the surface (and analytically depending on the point). Then the surface is an isotropic cyclide (defined in~\cite{krasauskas2014bilinear} and~\cite{sachs1990isotrope}).
\end{conjecture}

\begin{problem}
Determine the minimal number $k$ such that each surface in $\R^3$ containing $k$ isotropic circles through each point contains infinitely many ones.
\end{problem}

\begin{problem}
Can a surface parametrized by~\eqref{eqn:mainparam2} have degree exactly 5?
\end{problem}

The following conjecture gives a classification of surfaces containing a line and a parabola through each point.
\begin{conjecture}
Assume that through each point of an analytic surface in $\R^3$ one can draw a line segment and a parabolic arc with vertical axis fully contained in the surface, intersecting transversally, and depending analytically on the point. Then the surface (possibly besides a one-dimensional subset) has a parametrization
$$
\Phi(u,v)=\left(\frac{P}{UV},\frac{Q}{UV},\frac{Z}{U^2V}\right)
$$
for some $P,Q\in\R_{1,1},U\in\R_{1,0},V\in\R_{0,1},Z\in\R_{2,1}$, where $U,V\ne0$, such that the line segments are the curves $u=\const$ and the parabolic arcs are the curves $v=\const$.
\end{conjecture}

The following conjecture gives a complete description of the top views for the surfaces containing two isotropic circles through each point similar to Corollary~\ref{cor:top1}.
\begin{conjecture}\label{con:topview}
For each surface satisfying the assumptions of Conjecture~\ref{con:main2} the top views of the two isotropic circles through each point are either tangent to one cyclic (probably without real points) at two points (which may coincide or be complex conjugate) or lie in a union of two pencils of circles.
\end{conjecture}

\begin{example}[R.\,Krasauskas, private communication]
The surface given by the equation $x^2+y^2-(z-(x^2+y^2))^2=0$ contains \emph{exactly} two isotropic circles through each point. The top views of all isotropic circles on this surface form two pencils of circles.
\end{example}

Actually for a general position cyclic there are four families of circles which are tangent to this cyclic at two points (which may coincide or be complex conjugate). The other two families need not be top views of isotropic circles on the surface (R.\,Krasauskas, private communication)

Each parabola has a vertex and each isotropic circle that is not contained in a horizontal plane has points with minimal and maximal $z$-coordinates (called \emph{the lowest} and \emph{the highest} points respectively). Note that isotropic circles look like self-supporting arcs near these points (cf.~\cite{pottmann2007discrete,vouga2012design}). This is the reason why the following problem can be interesting for architecture.
\begin{problem}
For a surface containing two parabolas with vertical axes through each point describe the locus of parabolas' vertices. For a surface containing two isotropic circles through each point describe the loci of the lowest and the highest points of isotropic circles through each point.
\end{problem}

Another topic of interest for applications is rational B\'ezier representations.
\begin{problem}
Describe the control nets of biquadratic rational B\'ezier representations of surfaces containing two parabolas with vertical axes or two general isotropic circles through each point.
\end{problem}

Using isotropic model of Laguerre geometry, one can apply the results of the present paper to the following problem.
\begin{problem}\label{prb:laguerre}
Describe all surfaces such that there are two 1-parametric families of cones of revolution touching the surface along curves.
\end{problem}

\section*{Acknowledgments}
This work has been presented at the conferences ``Topology, Geometry, and Dynamics: Rokhlin --- 100'' in St.~Petersburg, ``Department of Higher Algebra becomes 90'' in Moscow, and 23rd M\"obius Contest. I am grateful to my scientific supervisor M.\,B.\,Skopenkov for stating the problem and constant attention to this work. I am also grateful to R.\,Krasauskas, N.\,Lubbes, and H.\,Pottmann for useful comments. Special thanks to A.\,A.\,Zaslavsky for his proof of Proposition~\ref{prp:two-circles}.


\end{document}